\documentclass[12pt]{article}

\usepackage{amsmath,amssymb,amsfonts} 

\usepackage{natbib} 
\usepackage{abstract} 
\usepackage{authblk} 
\usepackage[utf8]{inputenc} 

\usepackage{geometry}
\usepackage{graphicx} 
\usepackage[labelfont=bf]{caption}
\usepackage{subcaption}
\usepackage{algorithm}
\usepackage{algpseudocode}
\usepackage{hyperref}  
\usepackage{multirow} 
\usepackage{multicol} 
\usepackage{hhline}

\usepackage{tikz}
\usetikzlibrary{decorations.pathreplacing,calc}
\newcommand{\tikzmark}[1]{\tikz[overlay,remember picture] \node (#1) {};}
\newcommand*{\AddNote}[4]{%
    \begin{tikzpicture}[overlay, remember picture]
        \draw [decoration={brace,amplitude=0.5em},decorate,ultra thick,black]
            ($(#3)!(#1.north)!($(#3)-(0,1)$)$) --  
            ($(#3)!(#2.south)!($(#3)-(0,1)$)$)
                node [align=center, text width=5cm, pos=0.5, anchor=west] {#4};
    \end{tikzpicture}
}%
\begin{document}

\title{\textbf{Machine-learning-based arc selection for constrained shortest path problems in column generation}}
\author[1,2,3]{Mouad Morabit}
\author[1,2]{Guy Desaulniers}
\author[1,2,3,4]{Andrea Lodi}

\affil[1]{Department of Mathematics and Industrial
Engineering, Polytechnique Montréal (Québec)
Canada, H3C 3A7}
\affil[2]{GERAD, Montréal (Québec), Canada, H3T 2A7}
\affil[3]{Canada Excellence Research Chair in Data Science for Real-Time Decision-Making, Polytechnique Montréal}
\affil[4]{Jacobs Technion - Cornell Institute\\
Cornell Tech and Technion - IIT}
\affil[ ]{\textit {mouad.morabit@gerad.ca}}
\affil[ ]{\textit {guy.desaulniers@gerad.ca}}
\affil[ ]{\textit {andrea.lodi@cornell.edu}}
\date{\vspace{-5ex}}
\maketitle
    \begin{abstract}
Column generation is an iterative method used to solve a variety of optimization problems. It decomposes the problem into two parts: a master problem, and one or more pricing problems (PP). The total computing time taken by the method is divided between these two parts. In routing or scheduling applications, the problems are mostly defined on a network, and the PP is usually an NP-hard shortest path problem with resource constraints. In this work, we propose a new heuristic pricing algorithm based on machine learning. By taking advantage of the data collected during previous executions, the objective is to reduce the size of the network and accelerate the PP, keeping only the arcs that have a high chance to be part of the linear relaxation solution. The method has been applied to two specific problems: the vehicle and crew scheduling problem in public transit and the vehicle routing problem with time windows. Reductions in computational time of up to 40\% can be obtained.
\end{abstract}

\section{\label{sec:intro}Introduction}
Column generation (CG, \cite{desaulniers2005a}) is an iterative method that splits the original problem into two parts. The first part is a restricted master problem (RMP) that corresponds to the original linear problem (LP) but restricted to a subset of variables, which is often solved using the Simplex method. The second part is the subproblem, also called the pricing problem (PP). The role of the PP is to generate new columns of negative reduced cost (in case of a minimization problem) that can improve the current RMP solution. The method stops when no more columns can be generated. Otherwise, the new columns are added to the RMP to start a new iteration. To obtain integer solutions, CG is often embedded within a branch-and-bound process, where CG is used to solve the linear relaxation at each node of the tree. In this case the method is referred to as branch-and-price, see \cite{barnhart1998}, \cite{desaulniers2005a}.

Column generation has been successfully used to solve a variety of optimization problems, notably crew scheduling and vehicle routing problems arising in freight, urban and airline tranportation. For most of these problems, the master problem corresponds to a set partitioning problem with side constraints, and the PP is most likely a shortest path problem with resource constraints (SPPRC) or one of its variants. The success of using CG for this type of problems is also due to the efficient methods for solving the SPPRC, which are mainly dynamic programming algorithms \citep[]{desrochers-spprc-1988,irnich-spprc}. The goal is to find a path with the least cost between a source and a destination, while respecting the constraints on the resources (i.e., time, load of a vehicle, break duration, ...). It is well known that the SPPRC is NP-hard, but the methods developed to solve it are quite efficient (i.e., pseudopolynomial time complexity), and have allowed to solve large instances in reasonable times. The algorithms based on dynamic programming are also flexible enough to tackle different variants of the problems and take into account complex constraints that influence the feasibility of a route or a schedule. An improvement over the basic SPPRC algorithm is that of \citet{righini-bidir-espprc} who proposed a bi-directional search that propagates labels in both directions at the same time, i.e., from the source node to the destination, and backward from the destination to the source. Forward and backward labels residing in the same intermediate node are possibly merged to form complete feasible paths. 

The standard SPPRC often arises in problems that are defined on acyclic time-space networks, among others, in vehicle and crew scheduling problems \citep[]{Desaulniers1998,haase-desaulniers-desrosiers} and rostering problems \citep[]{gamache-rostering-99}. Another well-known variant that is more difficult to solve is the elementary SPPRC (ESPPRC) occuring in vehicle routing problems (VRP) \citep[]{Feillet2004, Contardo2015, costa-contardo-desaulniers-2019}. The use of an exact method to solve the ESPPRC is known to yield a tight lower bound, but given the difficulty of the problem, other alternatives based on relaxations and heuristics have been explored; see for example: ng-routes \citep{baldacci-ng-routes-2011}, SPPRC-k-cyc \citep{irnich-spprc-k-cyc} and partial elementarity \citep{desaulniers-lessard-tabu-kpath}. The goal of these later approaches is to find a trade-off between the quality of the bound and the difficulty of solving the PP. Several SPPRC heuristics have also been proposed in the literature. For instance, some authors relax the dominance rule applied when solving the problem with a labeling algorithm by only dominating on a subset of the resources \citep{desaulniers-lessard-tabu-kpath}. Another option is to keep only a limited number of labels at each node \citep{fukasawa2006}. Other proposed strategies are based on reducing the size of the network: one consists of sorting the incoming and outgoing arcs of each node by their reduced cost \citep{desaulniers-lessard-tabu-kpath}, where only a subset of the arcs with the least costs are retained. In \citet{fukasawa2006}, the authors use an extension of Kruskal's algorithm to build disjoint spanning trees from the network, and then, they consider only the arcs used on those trees in addition to the depot outgoing and incoming arcs. All of the above strategies can be used in the same execution and some can even be combined and used simultaneously. In most cases, an exact algorithm is invoked at the last iteration to prove the solution optimality. 

Combining machine learning methods with operations research and combinatorial optimization algorithms \citep{bengio2018} has been the subject of many studies over the last few years. These studies generally fall into two categories.
The first category is based on what is called \textit{imitation learning}, where the learner tries to imitate an expert in order to get a fast prediction to make decisions that are otherwise computationally expensive. In this category one seeks to reproduce decisions that are as close as possible to those of the expert. The presence of data from previous executions of the expert is therefore necessary.
The second category is \textit{reinforcement learning}, also called learning by experience as the name indicates, where the learner (i.e., agent) interacts with the environment and tries to optimize the results (i.e., maximize the rewards) of the decisions (i.e., actions) it makes.
In this context, the agent learns by trial and error while optimizing its reward function. After a number of iterations (i.e., epochs), a policy is learned telling the agent what is the best decision to make in each situation (i.e., state). The approach described in this article falls in the first category, that is learning by imitation and is thus based on supervised learning.

Several papers have covered the branch-and-bound method in mixed integer programming context, many of which seek to imitate \textit{strong branching} which is known to be computationally heavy but very effective, such as the works of \citet{Khalil2016, alvarez, gasse2019}. Several other methods exist and are covered in the survey of \citet{lodi2017}. For a more general overview of ML methods that have been developed and used in the context of combinatorial optimization (CO), readers are invited to check the survey of \citet{bengio2018}. This survey covers various machine learning methods and describes how they can be integrated to help solve CO problems. In the context of CG and branch-and-price, there are some very interesting papers worth mentioning. The authors in \citet{vaclavik2018} proposed a regression model for estimating an upper bound on the optimal value of the PP, which is then used to speed up the PP solution process. At each iteration of the CG, features based on the dual values are fed to the model. The learning is done in a so-called ``online" way (i.e., the learning is performed at the same time as the optimization), the loss function used during training penalizes the model according to the difference between the predicted bound and the true bound value (a high penalty is applied when the value of the bound is underestimated) and a discount factor is also used to give more importance to the last iterations. Also in the CG context, this time on the master problem side, a recent paper by \citet{mouad-columnselection} presents a learned column selector that selects, from a large set of generated columns, the columns to be added to the RMP at each CG iteration. The goal is to add the most promising columns that have a high probability of improving the current solution. The model used is based on graph neural networks, i.e., neural networks applied to graph-structured data. The problem is represented in a bipartite graph to model the variable-constraint relations. Differently from the current paper, the strategy in \citet{mouad-columnselection} is designed for problems that take a larger portion of the computing time in solving the RMP. 

As previously mentioned, CG can be applied to several optimization problems, so both RMP and PP can be different depending on the application at hand. In this article, the focus is more on crew scheduling and vehicle routing problems, where in most cases, the PP is a shortest path problem with resource constraints or one of its variants, and it is most likely the part that consumes the largest portion of the computing time. The purpose of this paper is to propose a heuristic pricing algorithm based on machine learning. The method takes advantage of the data collected from previous executions in order to learn a classifier that will select the arcs that have a greater chance to be used or to be part of an optimal solution. A good reduction of the network size decreases the computing time of the PP without increasing significantly the number of iterations, and thus reducing the overall total time. To ensure the exactness of the overall algorithm, the complete network is used when the reduced network fails to generate new columns. 

For our tests, we chose to use two well-known problems in the literature: the vehicle and crew scheduling problem (VCSP) and the vehicle routing problem with time windows (VRPTW). Each of these two problems have a different network structure (i.e., different type of nodes and arcs, etc), but nothing prevents from restricting the selection strategy to a specific type of arcs, making the methods described here 
general enough to be applied to various problems.

The remainder of the paper is organized as follows. Section 2 is devoted to some preliminaries and details about CG and the dynamic programming method used to solve the SPPRC. In Section 3, we describe the proposed selection strategy. Sections 4 and 5 will cover our two case studies, the VCSP and VRPTW respectively, highlighting the implementation differences between the two problems and reporting the computational results for each of them. Finally, we draw some conclusions in Section 6.

\section{\label{sec:prelim}Preliminaries}

In the CG context, and by focusing on the linear relaxation of the problem, let us consider the following master problem (MP) formulated as follows:

\begin{align}
z^{*}_{MP} := \min_x\hspace*{3mm} \sum_{p\in \Omega} c_px_p & \label{mp} \\
\mbox{s.t.} \hspace*{3mm} \label{const-partition1}
\sum_{p\in \Omega} \mathbf{a_p}x_p = \mathbf{b}&,  \\ \label{const-add1}
x_p\ge 0&,\hspace*{3mm} \forall p \in \Omega,
\end{align}

where $\Omega$ represents the set of variable indices (e.g., feasible routes or schedules), $c_p \in \mathbb{R}$ the variable cost, $\mathbf{a_p} \in \mathbb{R}^m$ the constraints coefficients and $\mathbf{b}\in \mathbb{R}^m$ the constraints right-hand side vector. When $|\Omega|$ is large and the variables cannot be enumerated explicitly, we consider only a subset $\mathcal{F} \subseteq \Omega$ of the variables, obtaining a restricted version of the problem above called a restricted master problem, RMP. At each CG iteration, the RMP is optimized, and an optimal solution $x$ is obtained along with a dual solution $\pi \in \mathbb{R}^m$ associated with the constraints~(\ref{const-partition1}). The dual values $\pi$ are then used to define the PP $\min_{p\in \Omega}\thickspace\{c_p - \pi^T\mathbf{a_p}\}$ and find new columns with negative reduced cost by solving it. The method stops when no such columns exist. In our applications, the variables represent either routes or schedules, and they can be generated by solving one or more PP (in our case a SPPRC or an ESPPRC).
\subsection{\label{subsec:spprc-form}SPPRC formulation}

Let $G=(V,A)$ be a directed graph with a set of nodes $V$ that include the source and destination nodes denoted by $s$ and $t$ respectively, and $A$ the set of arcs. Let $R$ be the set of resources. For each node $i \in V$, we define $T^r_i$ as the value of resource $r \in R$ accumulated on a partial path from source node $s$ to node $i$, and resource windows $[\underline{w}^r_i, \bar{w}^r_i]$, $\underline{w}^r_i, \bar{w}^r_i \in \mathbb{R}$ restricting the values that resource $r \in R$ can take on node $i$. For each arc $(i,j) \in A$, we define the resource consumptions $t_{ij}^r \in \mathbb{R}, r \in R$ and the cost $c_{ij}$. Let $X_{ij}$ be the decision variables of the model, which represent the flow on the arcs $(i,j)\in A$. The SPPRC is then formulated as follows:

\begin{align}
\min \hspace*{3mm} \sum_{(i,j)\in A} c_{ij}X_{ij} & \label{spprc_obj} \\
\mbox{s.t.} \hspace*{3mm} \label{flow_conservation_cst}
\sum_{j\in V} X_{ij} - \sum_{j\in V} X_{ji} =&
\begin{cases}
      +1 & \text{if $i = s$}\\
      0 & \text{if $i \in V$} \symbol{92} \{s,t\}\\
      -1 & \text{if $i = t$}
\end{cases} \\
\label{res_consump_cst}
X_{ij}(T_i^r + t_{ij}^r - T_j^r) \leq 0 &,\hspace*{3mm} \forall r \in R, \forall(i,j) \in A, \\
\label{res_windows_cst}
\underline{w}^r_i \leq T_i^r \leq \bar{w}^r_i &, \hspace*{3mm} \forall r \in R, \forall i \in V, \\
\label{binary_integer_cst}
X_{ij} \in \{0,1\}&, \hspace*{3mm} \forall (i,j) \in A,
\end{align}

where the objective~(\ref{spprc_obj}) minimizes the total cost of the path, the constraints~(\ref{flow_conservation_cst}) ensure the flow conservation along the path, while the constraints~(\ref{res_consump_cst}) model the resource consumption for each resource $r \in R$ on arc $(i,j) \in A$. Constraints~(\ref{res_windows_cst}) make sure that the resource values accumulated respect the resource intervals at each node and constraints~(\ref{binary_integer_cst}) are the binary requirements on the flow variables $X_{ij}$. Note that there are more complex SPPRC than those expressed by (\ref{spprc_obj})-(\ref{binary_integer_cst}), i.e., using complex resource extension functions that, given the resource values $T^r_i$ at node $i$, provide a lower bound on the resource values at node $j$.\\
As mentioned above, at each CG iteration, the dual values are used find new columns to add in the RMP. This dual information is included in the SPPRC by using a modified cost $\bar{c}_{ij} = c_{ij} - \sum_{k=1}^m \rho_{ij}^k \pi_k$ on each arc $(i,j) \in A$, where $\rho_{ij}^k$ is a binary parameter, taking value $1$ if the arc $(i,j) \in A$ contributes to the master problem constraint $k$, 0 otherwise. The cost of a path $p \in \Omega$ can be written in terms of the individual arcs:
\begin{align}
c_p = \sum_{(i,j) \in A} c_{ij}b_{ij}^p,
\end{align}
where $b_{ij}^p$ is a binary parameter equal to 1 if the arc $(i,j) \in A$ is included in path $p$. The constraints coefficients $a_p^k, k =\{1,2,\dots,m\}$, where $\mathbf{a_p} = [a_p^1, a_p^2, \dots, a_p^m]^T$ can also be written as
\begin{align}
a_p^k = \sum_{(i,j) \in A} \rho_{ij}^kb_{ij}^p.
\end{align}
The reduced cost $c_p$ of a variable/path can then be defined as the sum of the modified costs of the arcs composing the path, namely
\begin{align}
\bar{c}_p = c_p - \pi^T\mathbf{a_p} = c_p - \sum_{k=1}^m \pi_k a_p^k = \sum_{(i,j) \in A} b_{ij}^p(c_{ij} - \sum_{k=1}^m \rho_{ij}^k \pi_k) = \sum_{(i,j) \in A} b_{ij}^p\bar{c}_{ij}.
\end{align}

Therefore, solving the program (\ref{spprc_obj})-(\ref{binary_integer_cst}) using the modified costs in (\ref{spprc_obj}) yields a column with the least reduced cost. However, by solving the program using the dynamic programming formulation, it is possible to obtain several paths at once, which is another strong point for using these methods, as it has proven to speed up the resolution and minimize the number of iterations performed.

\subsection{\label{subsec:dyn-label-algo}Labeling algorithm}

In dynamic programming, shortest path problems are solved by a labeling algorithm. Such an algorithm starts from the trivial path that contains only the source node $s$, then extends it recursively along all arcs that yields a feasible path until reaching the destination node $t$. A path $P = (v_0, v_1, \dots, v_p)$ in $G$ represents a sequence of visited nodes such that $(v_{i-1},v_i) \in A, v_i \in V, \forall i=1,2, \dots, p$ where $v_p$ is the last node of the partial path and $v_0 = s$. In the algorithm, a label $L$ encodes information about a partial path $P$ such as the reduced cost $\bar{c}(L)$, the set of nodes $V(L)$ visited by the path, the last node visited by the path $v(L)$ and the accumulated quantity $T^r(L)$ of each resource $r \in R$. For a label $L$ ending at node $v_i (v_i \neq t)$, a new label $L'$ is obtained when an extension is performed along an arc $(i,j) \in A$, such that
\begin{align}
&v(L')=j, \\
&V(L')=V(L)\cup \{j\}, \\
&\bar{c}(L')= \bar{c}(L) + \bar{c}_{ij}, \\
&T^r(L') = \max\{ \underline{w}^r_j, T^r(L)+t_{ij}^r\}, \forall r \in R.
\end{align}

\begin{figure}[t]
\centering
\includegraphics{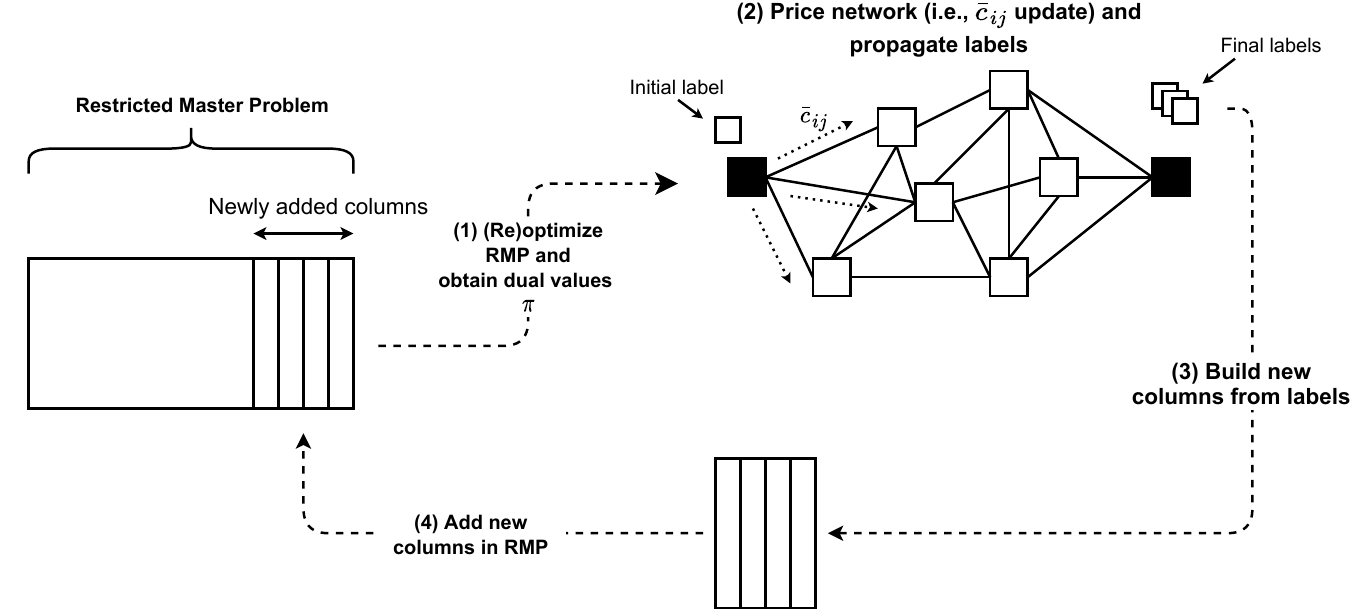}
    \caption{Summary of the different steps of a CG iteration.}
    \label{fig:one-cg-iteration}
\end{figure}

The label $L'$ also keeps a reference to its predecessor label $L$, in order to build the complete path when the algorithm ends. The label is considered infeasible if $T^r(L') > \bar w^r_j$ and is therefore deleted.\\
In addition, to the infeasible labels that get rejected when extending the labels, a dominance rule can be applied to eliminate non-promising paths. A label $L_1$ is said to be dominating a label $L_2$ if the following conditions are verified:
\begin{align}
v(L_1) = v(L_2),&\\
V(L_2) \subseteq V(L_1),&\\
\bar{c}(L_1) \leq \bar{c}(L_2),&\\
T^r(L_1) \leq T^r(L_2),& \hspace*{3mm} \forall r \in R
\end{align}
The performance of the labeling algorithm is closely linked to the choice of the dominance rule, a good choice of the latter allows to reduce the search space and to eliminate a maximum of labels, allowing to accelerate the resolution. Note that the validity of a dominance rule depends closely on the variant of the SPPRC being solved, as well as the resource consumptions on the arcs. For simplicity, some details have been omitted, such as resource extension functions (REFs, \citet{irnich2008}), which are considered being a more general way to define the resource arc consumptions, and extensions of the basic algorithm for the ESPPRC case to deal with cycle elimination. Readers are referred to \citet{irnich-spprc} for a good overview over the details omitted here.\\
At the end of the algorithm, several labels are obtained at the destination node $t$. The labels with a non-negative reduced cost are discarded, and those with a negative reduced cost are kept to build the new columns to be added in the RMP. The different steps discussed in this section are summarized in Figure \ref{fig:one-cg-iteration}.

\section{\label{sec:methodology}Methodology}

The goal of the project is to reduce the size of the underlying SPPRC network by keeping only the most promising arcs, therefore obtaining a reduced network. At each iteration, the same reduced network is used as long as it generates a satisfactory number of new negative reduced cost columns. The arc selection is thus performed only once before starting to solve the problem. In machine learning terms, this is a classification problem (i.e., supervised machine learning) and requires a labeled dataset 
$\mathcal{D} = \{(\boldsymbol{x}_1, y_1),(\boldsymbol{x}_2, y_2),\ldots,(\boldsymbol{x}_{n}, y_{n})\}$, where $n=|\mathcal{D}|$ and each entry is a tuple that represents an input/output pair, the inputs $\boldsymbol{x}_i$ are called features and $y_i$ are the desired outputs, also called labels (not the same notion of ``label" as in a labeling algorithm) given by an expert or known in advance. Given the dataset, the learning algorithm learns to predict the outputs $y_i$ from the feature vectors $\boldsymbol{x_i}$. In our case, it is a binary classification problem (i.e., $y_i \in \{0,1\}$), so each arc is either selected or not. During the training phase, the algorithm tries to minimize the difference between the true known values $y_i$ and the values predicted by the model denoted by $\hat{y_i}$. 
However, the goal is to be able to generalize and give good predictions for inputs never encountered during training, so a portion of the dataset is reserved for testing, in order to evaluate more accurately the performance of the model.
\subsection{\label{sec:features}Features}

The extracted features represent the characteristics of each individual arc $(i,j) \in A$. They can be different depending on the problem in hand and the structure of the underlying network. For the VCSP and VRPTW considered in this paper, some of the features are similar, namely:
\begin{itemize}
  \item Cost of the arc $c_{ij}$
  \item For each resource $ r \in R$, the resource consumption $T_{ij}^r$ along the arc $(i,j)$
  \item Number of arcs leaving node $i$
  \item Number of arcs entering node $j$
  \item For each resource $ r \in R$, the minimum, maximum and average resource consumptions along the arcs leaving node $i$
  \item For each resource $ r \in R$, the minimum, maximum and average resource consumptions along the arcs entering node $j$
  \item For each resource $ r \in R$, the upper and lower resource bounds for nodes $i$ and $j$.
\end{itemize}
The other features specific to each problem will be described in their dedicated sections.
\subsection{\label{sec:labels}Labels}

\begin{figure}[t]
\centering
\includegraphics[width=1.0\textwidth]{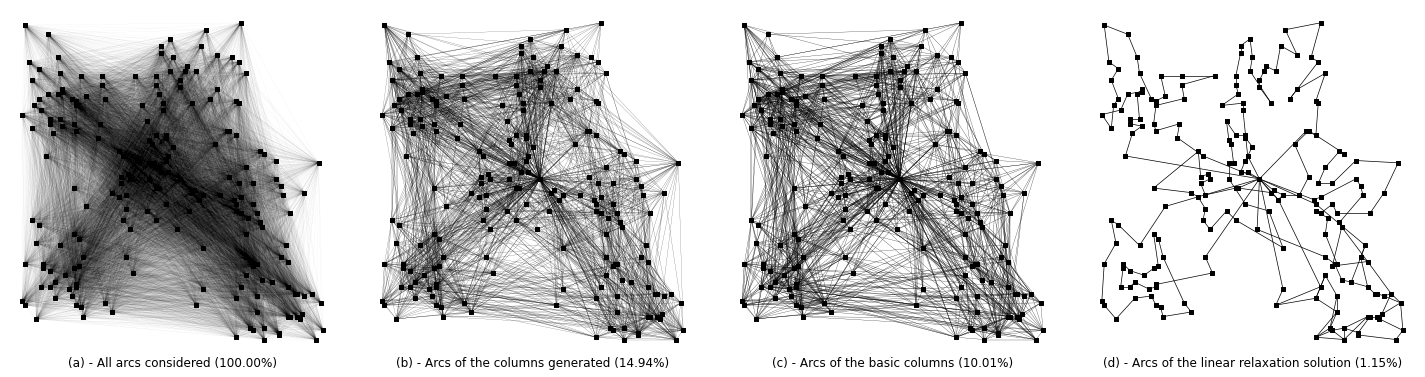}
    \caption{Number of arc comparisons across different scenarios}
    \label{fig:number_arcs_comparison}
\end{figure}

In supervised learning, the labels are either known in advance or assigned by an expert. In our case, we have to define what are the promising arcs to keep. The algorithm followed to assign the labels is in most cases computationally expensive, so we try to collect as much data as possible, then train a machine learning model that gives predictions in a reasonably fast computing time. Let $\mathcal{C}_i$ be the set of columns generated at iteration $i$, and $A_c$ the set of arcs visited by the path represented by column $c \in \mathcal{C}_i$. A promising arc is an arc that has been used at least once to generate columns at any iteration. The idea is that the generated columns must have dominated many other columns during the resolution of the PP and are therefore all good candidates, and so are the arcs composing their paths. According to the tests on the VCSP and VRPTW instances, the percentage of arcs used at least once during the resolution of the linear relaxation by CG does not exceed 18\%, and this number decreases when moving to larger instances.
Another idea is to consider only the positive basic columns (i.e., columns in the optimal basis of an RMP with a positive value), which results on a larger reduction of the graph, but according to the results, the number of iterations slightly increases compared to considering all the generated columns and the gain is therefore lost. Figure \ref{fig:number_arcs_comparison} shows a visual comparison of the number of arcs for a VRPTW instance, considering different scenarios: (a) all arcs, (b) arcs in the generated columns, (c) arcs in an RMP solution, and (d) arcs in the linear relaxation solution. The figure also specifies for each scenario the percentage of arcs selected with respect to the complete arc set.
\subsection{\label{sec:data_collection}Data collection}

Now that we have defined the arcs to be selected and the features, we can proceed and put all the pieces together to define our algorithm for data collection described in Algorithm \ref{alg:data_collection}.

\begin{algorithm}

\caption{Data collection : features and labels extraction}
\label{alg:data_collection}
\begin{algorithmic}[1] 

 \State $\mathcal{D} = \text{\O}$ 
 \label{stp1}
 \For{\texttt{each $a \in A$}} \hspace*{25mm} \label{stp2} \tikzmark{rightp1} \tikzmark{topp1}
 \State $y_a = 0$  
 \label{stp3} 
 \EndFor
 \label{stp4} \tikzmark{bottomp1}
 \State $\pi \longleftarrow$ \texttt{initialRMPSolution()}
 \label{stp4_2}
 \For{\texttt{each CG iteration $i$}} \label{stp5} 
     \State \texttt{priceNetwork(G,$\pi$)}
     \label{cg_1} 
     \State $\mathcal{C}_i \longleftarrow$ \texttt{generateColumns(G)}
     \label{cg_2} 
     \If{$\mathcal{C}_i == \text{\O}$} 
     \label{stop_1} 
     \hspace*{23mm} \tikzmark{rightp4} \tikzmark{topp4}
        \State \texttt{break}
        \label{stop_2}
     \EndIf
     \label{stop_3} \tikzmark{bottomp4}
     \For{\texttt{each $c \in \mathcal{C}_i$}} \hspace*{19mm} \label{stp6} 
     \tikzmark{rightp2} \tikzmark{topp2}
        \For{\texttt{each $a \in A_c$}} \label{stp7} 
        \State $y_a = 1$
        \label{stp8} 
        \EndFor
     \label{stp9} 
     \EndFor
     \label{stp10} \tikzmark{bottomp2}
     \State \texttt{addColumns($\mathcal{C}_i$)}
     \label{cg_3} 
     \State $\pi \longleftarrow$ \texttt{reoptimizeRMP()}
     \label{cg_4} 
 \EndFor
 \label{stp11} 
 \For{\texttt{each $a \in A$}} \hspace*{25mm} \label{stp12} \tikzmark{rightp3} \tikzmark{topp3}
    \State $\boldsymbol{x}_a \longleftarrow$ \texttt{extractData($a$)}
    \State $\mathcal{D} \longleftarrow \mathcal{D} \cup \{(\boldsymbol{x}_a, y_a)\}$ 
 \label{stp13} 
 \EndFor
 \label{stp14} \tikzmark{bottomp3}
 \end{algorithmic}
 \AddNote{topp1}{bottomp1}{rightp1}{Initialization}
 \AddNote{topp2}{bottomp2}{rightp2}{Label assignement to ``promising" arcs}
 \AddNote{topp3}{bottomp3}{rightp3}{Data collection}
 \AddNote{topp4}{bottomp4}{rightp4}{Solution is optimal}

\end{algorithm}

The steps of the algorithm are fairly straightforward. First, we start with an empty dataset in Step (\ref{stp1}) and we initialize all arc labels with the value 0 in Steps (\ref{stp2})-(\ref{stp4}). In (\ref{stp4_2}) an initial solution is obtained, either by using a heuristic or a two-phase simplex algorithm. Then, at each CG iteration, the RMP is reoptimized and the PP is solved, generating a new set of negative reduced cost columns $\mathcal{C}_i$ (Steps (\ref{cg_1})-(\ref{cg_2})). In case the PP defined over the full network G fails to generate any column, the CG process stops and the current solution is optimal (Steps (\ref{stop_1})-(\ref{stop_3})). Otherwise, for each newly generated column $c \in \mathcal{C}_i$, the label $y_a=1$ is assigned to each arc $a \in A_c$ composing the path represented by the column $c$ (Steps (\ref{stp6})-(\ref{stp10})) and the generated columns are added to the RMP which is repotimized (Steps (\ref{cg_3})-(\ref{cg_4})). Finally, when the CG process ends, for each arc $a \in A$ the pair of features $\boldsymbol{x}_a$ and labels $y_a$ are collected and added to the dataset $\mathcal{D}$ (Steps (\ref{stp12})-(\ref{stp14})). Note that it is necessary to solve several instances of the problem to optimality to collect enough data. The features can be extracted before solving the problem but for the labels it is necessary to wait until the end of the optimization. 

\subsection{\label{sec:ml_pricing}ML pricing algorithm}

Once the data is available, the preprocessing phase is performed (i.e., normalization of the values and the encoding of categorical features) before starting the training. Various classification algorithms can be used in our case: more details about the training phase and its results are provided in the next sections. Once the training phase is completed, the learned model can be used to select the arcs and build the reduced network. Before starting the first CG iteration, the learned model takes the features as an input and gives predictions $\hat{y}_a$ for each arc $a \in A$. The reduced network is defined as $G_r=(V, A_r)$ where $A_r = \{a \in A \hspace{1mm}|\hspace{1mm} \hat{y}_a = 1\}$ and it is used to generate columns as long as it yields a satisfactory number of columns, i.e., a number higher than a parameter value $\eta_{min}$. When the reduced network fails to generate enough columns, the complete network is used instead. If the latter generates a number of columns higher than a parameter value $\eta_{max}$, we switch back to the reduced network until we reach an optimal solution. Algorithm \ref{alg:cg_ml_pricing} details how the ML model is used in practice.

\begin{algorithm}

\caption{Machine-learning-based pricing heuristic in CG}
\label{alg:cg_ml_pricing}
\begin{algorithmic}[1] 
 \For{\texttt{each $a \in A$}} \hspace*{58mm} \label{alg2_stp1} \tikzmark{rightp1} \tikzmark{topp1}
 \State $\boldsymbol{x}_a \longleftarrow$ \texttt{extractFeatures($a$)} 
 \label{alg2_stp2} 
 \State $\hat{y}_a \longleftarrow $\texttt{predict($\boldsymbol{x}_a$)}
 \label{alg2_stp3}
 \EndFor
 \label{alg2_stp4}
 \State $A_r := \{a \in A \hspace{1mm}|\hspace{1mm} \hat{y}_a = 1\}$
 \label{alg2_stp5}
 \State $G_r := (V, A_r)$
 \label{alg2_stp6}
 \State $G_{active} \longleftarrow G_r$
 \label{alg2_stp7} 
 \State \texttt{useReducedG} $\longleftarrow$ \texttt{True}
 \label{alg2_stp8} \tikzmark{bottomp1}
 \State $\pi \longleftarrow$ \texttt{initialRMPSolution()}
 \label{alg2_stp9}
 \For{\texttt{each CG iteration $i$}} \label{alg2_stp10} 
     \State \texttt{priceNetwork($G_{active}$,$\pi$)}
     \label{alg2_stp11}
     \State $\mathcal{C}_i \longleftarrow$ \texttt{generateColumns($G_{active}$)}
     \label{alg2_stp12}
     \If{($|\mathcal{C}_i| < \eta_{min}$ $\wedge$ \texttt{useReducedG})} \hspace*{20mm} \tikzmark{rightp2} \tikzmark{topp2}
     \label{alg2_stp13} 
        \State \texttt{useReducedG} $\longleftarrow$  $\texttt{False}$
        \label{alg2_stp14}
        \State $G_{active} \longleftarrow G$
        \label{alg2_stp15}

     \ElsIf{$|\mathcal{C}_i| \geq \eta_{max}$ $\wedge$ ($\lnot\texttt{useReducedG}$)}
     \label{alg2_stp16}
        \State \texttt{useReducedG} $\longleftarrow$  $\texttt{True}$
        \label{alg2_stp17}
        \State $G_{active} \longleftarrow G_r$
        \label{alg2_stp18}
     \tikzmark{bottomp2}
     
     \ElsIf{$\mathcal{C}_i == \text{\O}$ $\wedge$ ($\lnot\texttt{useReducedG}$)} \hspace*{12mm} \tikzmark{rightp3} \tikzmark{topp3}
     \label{alg2_stp19}
        \State \texttt{Exit}
        \label{alg2_stp20}
     \EndIf
     \label{alg2_stp21}
     \tikzmark{bottomp3}
     
     \If{$|\mathcal{C}_i| > 0$}
     \label{alg2_stp22}
         \State \texttt{addColumns($\mathcal{C}_i$)}
         \label{alg2_stp23}
         \State $\pi \longleftarrow$ \texttt{reoptimizeRMP()}
         \label{alg2_stp24}
     \EndIf
     \label{alg2_stp25}
 \EndFor
 \label{alg2_stp26}
 
 \end{algorithmic}
 \AddNote{topp1}{bottomp1}{rightp1}{Reduced network $G_r$}
 \AddNote{topp2}{bottomp2}{rightp2}{Network switch cases}
 \AddNote{topp3}{bottomp3}{rightp3}{Solution is optimal}
 
\end{algorithm}

The algorithm starts by extracting the features of the arcs and getting the model predictions (Steps (\ref{alg2_stp1})-(\ref{alg2_stp4})). The predictions $\hat{y}_a$ are then used to build the reduced graph $G_r$ that is set as the active graph to be used when solving the PP (Steps (\ref{alg2_stp5})-(\ref{alg2_stp8})). Note that the active graph $G_{active}$ can either be the original complete graph $G$ or the reduced one $G_r$. The boolean variable $useReducedG$ is also initialized to $\texttt{True}$, indicating that the reduced graph is being used (e.g., equivalent to $G_{active}==G_r$). In Step (\ref{alg2_stp9}), an initial solution to the RMP is obtained along with dual values $\pi$. At each CG iteration, the currently active network is priced and a new set of columns $\mathcal{C}_i$ is generated (Steps (\ref{alg2_stp11}) and (\ref{alg2_stp12})). If the columns were generated by the reduced network and the number of columns $|\mathcal{C}_i|<\eta_{min}$, we switch to using the full network by setting the variables $G_{active}$ to $G$ and $useReducedG$ to $\texttt{False}$ (Steps (\ref{alg2_stp13})-(\ref{alg2_stp15})). On the contrary, if the full network was used and $|\mathcal{C}_i|\geq\eta_{max}$, we switch back to using the reduced network $G_r$ and set $usedReducedG$ to $\texttt{True}$ (Steps (\ref{alg2_stp16})-(\ref{alg2_stp18})).
By using the two previous conditions, the algorithm switches between the two graphs until an optimal solution is obtained. Normally, with the right choice of the two parameters $\eta_{min}$ and $\eta_{max}$, at the last iterations of the CG, the reduced network will not generate enough columns (i.e., less than $\eta_{min}$) and so the full network $G$ will be used instead. Most likely, the full network will neither generate many columns (i.e., less than $\eta_{max}$), so the algorithm will not switch back to the reduced network and will keep using the full network until it generates no more columns, indicating that the obtained solution is optimal (Steps (\ref{alg2_stp19})-(\ref{alg2_stp21})). If columns were generated, they are added to the RMP which is reoptimized afterwards (Steps (\ref{alg2_stp22})-(\ref{alg2_stp25})).

\section{\label{sec:case1}Application I: Vehicle and crew scheduling problem}

In a public transit system, the planning process goes through several stages. The first step in the process is to determine the bus lines, the stops on each line, and the frequency of the trips. Based on these pieces of information, a timetable is created, describing the trips with their corresponding starting and ending times and locations. The next two steps in the process are the construction of vehicle routes and crew schedules, which means solving, respectively, the two scheduling problems: the \textit{Vehicle scheduling problem} and the \textit{Crew scheduling problem}. Traditionally, these two problems are solved in a sequential manner, where the vehicle routes are determined first before the crew schedules. However, this approach is not guaranteed to provide the best solution, because in most cases driver costs dominate the costs of vehicle use. A first formulation that integrates and considers the two problems simultaneously was proposed by \citet{freling1999}, giving rise to the VCSP. For solving the problem, the authors described a column generation approach applied to a Lagrangian relaxation of the master problem, since the problem contains a very large number of variables representing the feasible duties. Most of the formulations proposed in the literature are set partitioning ones, where column generation plays an important role to deal with the large number of variables and to generate schedules as needed. Given the difficulty of the problem, several heuristics have been proposed in the literature as well. In this paper, we consider the formulation described in \citet{haase-desaulniers-desrosiers}, where the RMP is the linear relaxation of a set partitioning problem with side constraints and the PP is a SPPRC. The RMP formulation contains only the crew schedule variables, but side constraints are added to the problem to ensure that the vehicle schedules can be obtained afterwards in a polynomial time. In addition, the costs of the vehicles are included in the objective function to ensure that an overall optimal solution is obtained. 

\subsection{\label{subsec:case1_network_structure}Network structure}

Since the SPPRC is the part of interest for our work, we will take a closer look at the network structure of the problem before diving into the details about the data collection, the instances used and the results obtained.

\begin{figure}[t]
\centering
\includegraphics[width=1.0\textwidth]{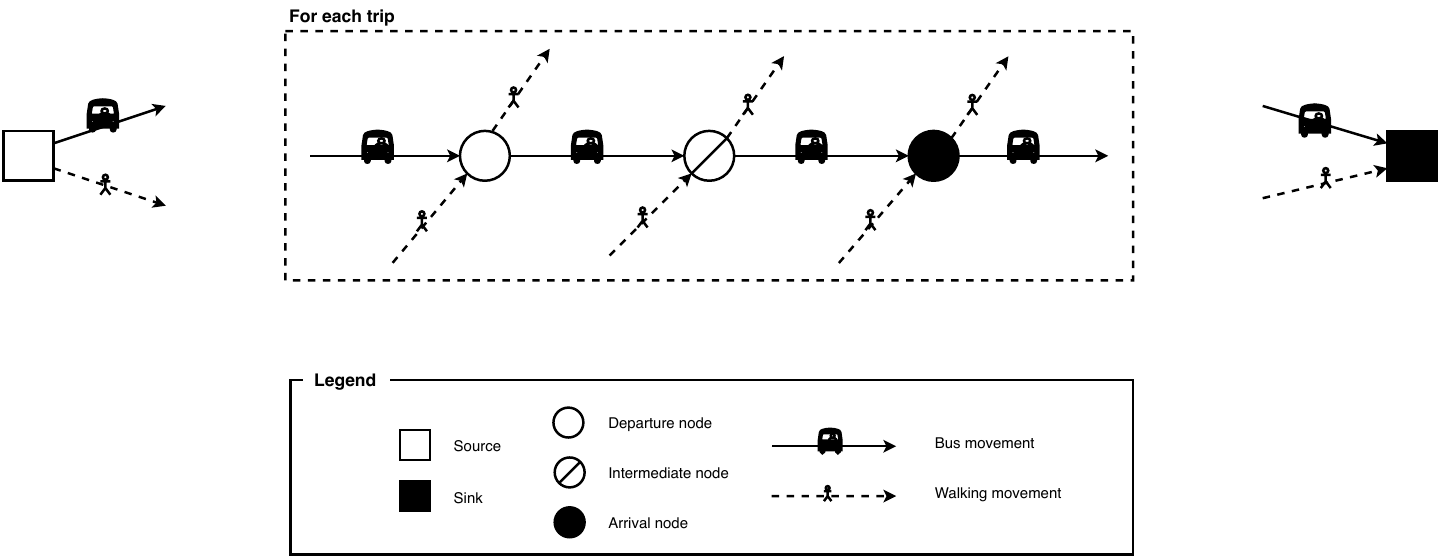}
    \caption{A simplified version of the network components in the VCSP.}
    \label{fig:network_structure}
\end{figure}

In this problem, a network is used to generate driver schedules and there is one network for each schedule type (e.g, regular schedule with one meal break or straight without a meal break). Additionally, the arcs can be divided into two categories: 1) bus movement when the driver is driving or attending a bus and 2) walking when the driver is moving on his/her own for repositioning purposes. As mentioned before, each bus line is defined by three different trip nodes: departure node, arrival node and intermediate nodes where an exchange of drivers can be performed. For each bus line, several trips are planned during the day at different times, the exact time at which a bus should arrive at each of the trip nodes is therefore known. To cover the trips, the buses leave the depot in the direction of the trip departure nodes. Each bus can be used for multiple trips, and can be operated by different drivers before returning to the depot at the end of the day. Since the same bus must service the whole trip, a bus moving to the departure node of a trip must be retained until it arrives at the destination node, at which point it can head to the start of the next trip. Therefore, no bus movements leaving the departure or intermediate nodes to a different trip are allowed. In fact, the bus movements represent less than 7\% of the total number of arcs, while the majority of the remaining 93\% are the walking movements of the drivers. Unlike the bus movements, the walking movements are not as restricted because they can occur at the beginning, the end or the middle of the trips, as long as they are feasible. Figure \ref{fig:network_structure} is a simplified and non-detailed representation of the main components found in a VCSP. In addition to the main components, the network can also include break nodes and their corresponding arcs to model the breaks that drivers can take after a certain working time. Some details are not involved in our selection strategy and therefore not detailed here, but the interested readers can consult \citet{haase-desaulniers-desrosiers} for an in-depth overview.

For us, the movements we are interested in are the walking movements of the drivers between trip nodes, since they represent more than 90\% of the total number of arcs. The selection strategy only covers this subset of arcs, denoted as $A_s \subset A$ and not all arcs in $A$. For a walking movement arc linking the trip nodes of two different trips $T1$ and $T2$, the classification model will try to predict whether there is a high chance that a driver leaves his/her current trip $T1$ in the direction of trip $T2$. The arcs entering and leaving the depot, as well as the bus movements entering and leaving the trip nodes are not involved, in addition to any extra arcs for breaks, etc.

\subsection{\label{subsec:case1_data_collection}Data collection, features and labels}

All the features mentioned in Section \ref{sec:features} are valid for the VCSP. For the resource consumptions on the arcs, the two representing the working time and waiting time, are respectively extracted. The first resource is the time taken to walk to the next trip node, while the second resource is the waiting time required before actually starting the trip. Note that the waiting time is bounded and cannot exceed the upper bound indicated by the resource windows. Given that the exact time at which the trip nodes must be visited is known, the times associated with node $i$ and $j$ for each arc $(i,j) \in A_s$ are added to the list of features as well. The difference between $j$ and $i$ times is equal to the time required to walk from $i$ to $j$, plus the waiting time, but since there are more trips during peak hours, the time of the day information can be useful. Another feature that can be added to the list is the arc type since the arcs we are interested in can link three different types of nodes (i.e., departure, intermediate and arrival nodes) yielding nine possible connection types. Therefore, with the additional features, a total of 22 features is collected for each arc.

As described in Section \ref{sec:labels}, labels are assigned to the arcs that are part of the generated columns. However, only the arcs representing the walking movements between the trip nodes are concerned.

\subsection{\label{subsec:case1_instances}VCSP instances}

To collect enough data, an instance generator is used as in \citet{haase-desaulniers-desrosiers}. By taking a total number of trips as an input, the generator splits this number over the different bus lines. The trips are then randomly distributed over the time horizon, with more trips during peak hours. A total of 20 instances with 400 trips were used in the training phase, each instance added about 250,000 data points to the dataset. Other new instances of different sizes (i.e., $300, 350, \dots, 600$) are generated and used to test the integration of the ML model in the CG algorithm. The instances of size other than $400$ also serve to evaluate the generalization of the model to unseen instances of different sizes.

\subsection{\label{subsec:case1_results}Computational results}

This section is divided into three parts. The first one is devoted to the machine learning phase, where additional details about data collection, the algorithms used and the results obtained are given, whereas the second part presents the results of the integration of the model in the CG algorithm. In the third part, we highlight the results of some additional experiments that are worth mentioning. All the experiments were performed on a Linux machine with an i7-8700 CPU @ 3.20GHz and 64GB of RAM.

\subsubsection{\label{subsubsec:case1_ml_results}Machine learning}

With the dataset in hand, some preprocessing is performed such as data normalization and one-hot encoding of categorical features (i.e., arc type). The next step is to choose the appropriate classification algorithm to use. There are several possibilities, from the simplest linear model (i.e., Logistic regression) to the more classic algorithms (SVM, Decision Trees, etc.) and more complex models (Artificial Neural Networks, etc). To compare the models with each other, the main criteria is their accuracy on the test set. However, when there is an imbalance between the output classes (i.e., in our case, only about 15\% of the data have label $1$ and the remaining have label $0$), it is better to consider other metrics such as the true positive rate (TPR, also called the Recall), the true negative rate (TNR), precision, etc. These metrics provide an idea about the performance obtained for each individual class allowing a better comparison between the different models.

Three different algorithms have been chosen for the training phase, starting with a linear classifier: Logistic regression, then moving to two more advanced algorithms: Random Forest (RF) and Artificial Neural Networks (ANN). Each of these algorithms requires hyperparameters to be tuned, a ``cross validation" approach is used to compare the different set of values. The best values obtained are described in Table \ref{tab:vcsp-ml-hyperparameter}, while Table \ref{tab:vcsp-metrics} shows the results of the three algorithms on the test set.

\begin{table}[t]
  \begin{center}
    
    \begin{tabular}{|c|l|c|} 
      \hline
      \textbf{Algorithm} & \textbf{Hyperparameter} & \textbf{Value} \\
      \hline
      \multirow{3}{*}{\textbf{Logistic regression}} & C parameter & 1\\
       & Solver & lbfgs\\
       & Class weights & Balanced\\
      
      \hline
      \multirow{7}{*}{\textbf{Random forest}} & Bootstrap & True\\
       & Max depth & 10\\
       & Max features & 7\\
       & Min samples per leaf & 50\\
       & Min samples per split & 100\\
       & Number of trees & 500\\
       & Class weights& Balanced\\
      
      \hline
      \multirow{8}{*}{\textbf{Neural Network}} & Learning rate & $10^{-3}$\\
       & Epochs & $1000$\\
       & Batch size & $32$\\
       & Loss function & Binary cross entropy\\
       & Activation function & ReLU\\
       & Architecture & 32x32x32x1\\
       & Optimizer & Adam\\
       & Class weights& 7:1\\
      \hline
    \end{tabular}
  \end{center}
  \caption{Hyperparameters values used in the training phase.}
  \label{tab:vcsp-ml-hyperparameter}
\end{table}

\begin{table}[t]
  \begin{center}
    
    \begin{tabular}{|l|c|c|c|} 
      \hline
      \textbf{Algorithm} & \textbf{Recall} & \textbf{TNR} & \textbf{Balanced accuracy} \\
      \hline
      Logistic regression & 64\% & 81\% & 72\%\\
      \hline
      Random forest & 76\% & 81\% & 78\%\\
      \hline
      Neural network & 78\% & 80\% & 79\%\\
      \hline
    \end{tabular}
  \end{center}
  \caption{VCSP metrics of three different algorithms.}
  \label{tab:vcsp-metrics}
\end{table}

The Recall represents the percentage of arcs with label $y_a=1$ that are correctly classified, while the TNR represents the accuracy of the opposite class (i.e., classification accuracy of the arcs with label $0$). Ideally, a high percentage of both is sought. The last column represents the balanced accuracy which is the average of the two previous columns. To deal with the unbalanced data, weights are used during the optimization of the loss function. The goal is to give a higher penalty to the arcs with label 1 that are misclassified since they are a minority in the dataset, hence the use of the ``balanced" setting in the class weights hyperparameter. One can imagine the dataset as a cloud of points, where 85\% are red and 15\% are green. With logistic regression, one tries to separate these two classes with a linear separator, but since there are more red points than green, there is more chance to have a large number of red points on one side of the separator and thus a high TNR. A high recall can also be obtained if the two classes are easily separable, which does not seem to be the case according to the results obtained (i.e., 64\% Recall and 81\% TNR). The model slightly underfits but overall, the accuracy obtained is not far from the other two algorithms. On the other hand, the RF and ANN models are able to reach a higher accuracy of 78-79\%, with more balance between the Recall and the TNR. The values obtained by the two algorithms are quite close to each other.

\subsubsection{\label{subsubsec:case1_cg_results}CG with arc selection}

After the training phase, both the RF and ANN models seem to be good candidates to be integrated in the CG algorithm. Given their almost similar accuracy, one can expect them to give comparable results. When testing the models on several groups of instances with different sizes, the RF model performs better on some instances, but on others it is the ANN model that is better. However, because the RF model is a few percents ahead on average, we retained this candidate for the subsequent tests. 

As described in Algorithm \ref{alg:cg_ml_pricing}, the model is used before the beginning of the CG process to build the reduced network $G_r$. Different scenarios are compared to assess the efficiency of the arc selection, namely:

\textbf{Baseline-CG}: This corresponds to the basic CG with no arc selection, i.e., network $G$ is used in all iterations.

\textbf{ML-S}: Arc selection is used and the CG alternates between the use of the reduced network $G_r$ and the full network $G$. The reduced network is built from the predictions obtained by the RF model. The parameters used during the CG process are: $\eta_{min}=30, \eta_{max}=100$.

\textbf{Random-S}: This scenario corresponds to a completely random selection of the arcs in $A_s$, the number of selected arcs is the same as with the model. This can be achieved by simply shuffling the values obtained by the model.

\textbf{Cost-S}: For each trip node $i$, incoming arcs $(n,i) \in A_s$ and outgoing arcs $(i,n) \in A_s$ are sorted in ascending order of their cost, and a subset of arcs with the lowest costs is selected. The total number of selected arcs is approximately the same as with the model.

The results obtained by the four algorithms are reported in Table \ref{tab:table-results-vcsp}. The names of the instances are written in the form ``VCS\_\textit{Size}\_\textit{Id}", where \textit{Size} is the instance size, i.e., the number of trips, and \textit{Id} is the instance identifier. For each algorithm, we report the number of CG iterations required to reach the linear relaxation solution, the computing time in seconds of: the PP, the RMP and in total. For the algorithms involving an arc selection, an additional column shows the reduction of the computing time gained in comparison with the ``Baseline-CG" algorithm. Additionally, the number of iterations performed using the full network $G$ is indicated in parentheses. After each group of instances of the same size, a line indicating the average values is added.

\begin{table}
  \begin{center}
    \resizebox{\textwidth}{!}{
    \begin{tabular}{|r||r|r|r|r||c|r|r|r|r||c|r|r|r|r||c|r|r|r|r|}
      \hline
      \multirow{3}{*}{\textbf{Instance}} &   \multicolumn{4}{|c||}{\textbf{Baseline-CG}} &  \multicolumn{5}{|c||}{\textbf{ML-S}} & 
      \multicolumn{5}{|c||}{\textbf{Random-S}} &
      \multicolumn{5}{|c|}{\textbf{Cost-S}} \\
      \cline{2-20}
      & \multirow{2}{*}{\#Itr} & \multicolumn{3}{|c||}{Time (s)}
      & \multirow{2}{*}{\#Itr} & \multicolumn{3}{|c|}{Time (s)} & \multirow{2}{*}{\textbf{Gain}} 
      & \multirow{2}{*}{\#Itr} & \multicolumn{3}{|c|}{Time (s)} & \multirow{2}{*}{\textbf{Gain}}
      & \multirow{2}{*}{\#Itr} & \multicolumn{3}{|c|}{Time (s)} & \multirow{2}{*}{\textbf{Gain}}
      \\
      \cline{3-5} \cline{7-9} \cline{12-14} \cline{17-19}
      & & PP & RMP & Total
      & & PP & RMP & Total & 
      & & PP & RMP & Total &
      & & PP & RMP & Total & \\
      
      \hline\hline
      \textbf{VCS\_300\_1} & 158 & 133 & 80 & 213
                           & 146 (12) & 75 & 76 & 151 & \textbf{29\%}
                           & 180 (59) & 114 & 91 & 205 & \textbf{4\%}
                           & 168 (34) & 94 & 82 & 176 & \textbf{17\%}\\
      \textbf{VCS\_300\_2} & 154 & 174 & 85 & 259
                           & 158 (8) & 106 & 91 & 197 & \textbf{24\%}
                           & 242 (72) & 214 & 122 & 336 & \textbf{-30\%}
                           & 183 (37) & 143 & 94 & 237 & \textbf{8\%}\\
      \textbf{VCS\_300\_3} & 157 & 182 & 81 & 263
                           & 154 (14) & 115 & 76 & 191 & \textbf{27\%}
                           & 204 (70) & 198 & 100 & 298 & \textbf{-13\%}
                           & 182 (37) & 160 & 86 & 246 & \textbf{6\%}\\
      \textbf{VCS\_300\_4} & 176 & 222 & 81 & 303
                           & 185 (12) & 166 & 81 & 247 & \textbf{18\%}
                           & 227 (105) & 220 & 96 & 316 & \textbf{-4\%}
                           & 220 (66) & 241 & 84 & 325 & \textbf{-7\%}\\
      \textbf{VCS\_300\_5} & 199 & 204 & 89 & 293
                           & 191 (27) & 133 & 87 & 220 & \textbf{25\%}
                           & 238 (60) & 184 & 95 & 279 & \textbf{5\%}
                           & 192 (44) & 166 & 87 & 253 & \textbf{14\%}\\
      \hline
      \textbf{Average} 
      & \textbf{169} & \textbf{183} & \textbf{83} & \textbf{266} 
      & \textbf{167 (14)} & \textbf{119} & \textbf{82} & \textbf{201} & \textbf{25\%} 
      & \textbf{218 (73)} & \textbf{186} & \textbf{101} & \textbf{287} & \textbf{-13\%} 
      & \textbf{189 (43)} & \textbf{161} & \textbf{87} & \textbf{248} & \textbf{8\%}\\
      \hline
      \hline
      \textbf{VCS\_350\_1} & 136 & 183 & 126 & 309
                           & 143 (15) & 112 & 129 & 241 & \textbf{22\%}
                           & 199 (67) & 207 & 151 & 359 & \textbf{-16\%}
                           & 175 (42) & 166 & 144 & 310 & \textbf{0\%}\\
      \textbf{VCS\_350\_2} & 148 & 158 & 136 & 294
                           & 135 (9) & 85 & 119 & 204 & \textbf{31\%}
                           & 184 (67) & 146 & 142 & 288 & \textbf{2\%}
                           & 192 (52) & 153 & 150 & 303 & \textbf{-3\%}\\
      \textbf{VCS\_350\_3} & 182 & 253 & 179 & 432
                           & 170 (10) & 141 & 167 & 308 & \textbf{29\%}
                           & 219 (61) & 226 & 207 & 433 & \textbf{0\%}
                           & 192 (37) & 192 & 178 & 370 & \textbf{14\%}\\
      \textbf{VCS\_350\_4} & 171 & 276 & 162 & 438
                           & 188 (7) & 173 & 174 & 347 & \textbf{21\%}
                           & 302 (97) & 362 & 243 & 605 & \textbf{-38\%}
                           & 227 (51) & 291 & 173 & 464 & \textbf{-6\%}\\
      \textbf{VCS\_350\_5} & 201 & 360 & 164 & 524
                           & 234 (12) & 272 & 183 & 455 & \textbf{13\%}
                           & 273 (101) & 392 & 222 & 614 & \textbf{-17\%}
                           & 249 (48) & 347 & 179 & 526 & \textbf{0\%}\\
      \hline
      \textbf{Average} 
      & \textbf{168} & \textbf{246} & \textbf{153} & \textbf{399} 
      & \textbf{174 (10)} & \textbf{157} & \textbf{154} & \textbf{311} & \textbf{23\%} 
      & \textbf{215 (78)} & \textbf{267} & \textbf{193} & \textbf{460} & \textbf{-5\%} 
      & \textbf{207 (46)} & \textbf{230} & \textbf{165} & \textbf{395} & \textbf{1\%}\\
      \hline
      \hline
      \textbf{VCS\_400\_1} & 177 & 370 & 253 & 623
                           & 177 (10) & 238 & 258 & 496 & \textbf{20\%}
                           & 212 (89) & 361 & 293 & 654 & \textbf{-5\%}
                           & 216 (38) & 342 & 297 & 639 & \textbf{-3\%}\\
      \textbf{VCS\_400\_2} & 176 & 422 & 246 & 668
                           & 189 (8) & 255 & 258 & 513 & \textbf{23\%}
                           & 240 (86) & 471 & 298 & 769 & \textbf{-15\%}
                           & 214 (58) & 378 & 257 & 635 & \textbf{5\%}\\
      \textbf{VCS\_400\_3} & 153 & 266 & 242 & 508
                           & 163 (9) & 177 & 233 & 410 & \textbf{19\%}
                           & 176 (64) & 263 & 238 & 501 & \textbf{1\%}
                           & 184 (42) & 245 & 254 & 499 & \textbf{2\%}\\
      \textbf{VCS\_400\_4} & 147 & 209 & 206 & 415
                           & 147 (12) & 126 & 190 & 316 & \textbf{24\%}
                           & 199 (64) & 225 & 224 & 449 & \textbf{-8\%}
                           & 173 (30) & 173 & 222 & 395 & \textbf{5\%}\\
      \textbf{VCS\_400\_5} & 207 & 518 & 256 & 774
                           & 180 (12) & 283 & 230 & 513 & \textbf{34\%}
                           & 266 (100) & 558 & 311 & 869 & \textbf{-12\%}
                           & 214 (44) & 391 & 252 & 643 & \textbf{17\%}\\
      \hline
      \textbf{Average} 
      & \textbf{172} & \textbf{357} & \textbf{241} & \textbf{598} 
      & \textbf{171 (10)} & \textbf{216} & \textbf{234} & \textbf{450} & \textbf{24\%} 
      & \textbf{218 (80)} & \textbf{376} & \textbf{273} & \textbf{648} & \textbf{-8\%} 
      & \textbf{200 (42)} & \textbf{306} & \textbf{256} & \textbf{562} & \textbf{5\%}\\
      \hline
      \hline
      \textbf{VCS\_450\_1} & 340 & 1114 & 640 & 1754
                           & 368 (14) & 809 & 636 & 1445 & \textbf{18\%}
                           & 426 (138) & 1034 & 713 & 1747 & \textbf{0\%}
                           & 366 (64) & 868 & 637 & 1505 & \textbf{14\%}\\
      \textbf{VCS\_450\_2} & 269 & 1540 & 544 & 2084
                           & 241 (11) & 797 & 486 & 1283 & \textbf{38\%}
                           & 358 (126) & 1453 & 663 & 2116 & \textbf{-2\%}
                           & 281 (45) & 1071 & 565 & 1636 & \textbf{21\%}\\
      \textbf{VCS\_450\_3} & 194 & 570 & 396 & 966
                           & 194 (10) & 319 & 391 & 710 & \textbf{27\%}
                           & 243 (81) & 497 & 463 & 960 & \textbf{1\%}
                           & 251 (64) & 508 & 435 & 943 & \textbf{2\%}\\
      \textbf{VCS\_450\_4} & 209 & 505 & 387 & 892
                           & 234 (8) & 312 & 440 & 752 & \textbf{16\%}
                           & 264 (110) & 473 & 446 & 919 & \textbf{-3\%}
                           & 255 (52) & 439 & 439 & 878 & \textbf{2\%}\\
      \textbf{VCS\_450\_5} & 279 & 1507 & 551 & 2058
                           & 289 (21) & 976 & 564 & 1540 & \textbf{25\%}
                           & 376 (124) & 1566 & 681 & 2247 & \textbf{-9\%}
                           & 350 (61) & 1510 & 568 & 2078 & \textbf{-1\%}\\
      \hline
      \textbf{Average} 
      & \textbf{258} & \textbf{1047} & \textbf{504} & \textbf{1551} 
      & \textbf{265 (12)} & \textbf{643} & \textbf{503} & \textbf{1146} & \textbf{25\%} 
      & \textbf{333 (115)} & \textbf{1005} & \textbf{593} & \textbf{1598} & \textbf{-3\%} 
      & \textbf{300 (57)} & \textbf{879} & \textbf{529} & \textbf{1408} & \textbf{8\%}\\
      \hline
      \hline
      \textbf{VCS\_500\_1} & 277 & 1618 & 856 & 2474
                           & 304 (8) & 1030 & 902 & 1932 & \textbf{22\%}
                           & 387 (132) & 1738 & 1021 & 2759 & \textbf{-12\%}
                           & 308 (60) & 1272 & 916 & 2188 & \textbf{12\%}\\
      \textbf{VCS\_500\_2} & 296 & 1654 & 620 & 2274
                           & 270 (13) & 846 & 576 & 1422 & \textbf{37\%}
                           & 366 (130) & 1404 & 738 & 2142 & \textbf{6\%}
                           & 342 (76) & 1318 & 691 & 2009 & \textbf{12\%}\\
      \textbf{VCS\_500\_3} & 203 & 770 & 604 & 1374
                           & 221 (10) & 453 & 611 & 1064 & \textbf{23\%}
                           & 291 (101) & 887 & 736 & 1623 & \textbf{-18\%}
                           & 266 (54) & 688 & 713 & 1401 & \textbf{-2\%}\\
      \textbf{VCS\_500\_4} & 266 & 1839 & 748 & 2587
                           & 281 (8) & 1235 & 749 & 1984 & \textbf{23\%}
                           & 357 (130) & 1842 & 908 & 2750 & \textbf{-6\%}
                           & 313 (54) & 1526 & 789 & 2315 & \textbf{11\%}\\
      \textbf{VCS\_500\_5} & 322 & 1720 & 627 & 2347
                           & 340 (23) & 1024 & 656 & 1680 & \textbf{28\%}
                           & 442 (173) & 1635 & 769 & 2404 & \textbf{-2\%}
                           & 386 (103) & 1486 & 701 & 2187 & \textbf{7\%}\\
      \hline
      \textbf{Average} 
      & \textbf{273} & \textbf{1520} & \textbf{691} & \textbf{2211} 
      & \textbf{283 (12)} & \textbf{918} & \textbf{699} & \textbf{1616} & \textbf{27\%} 
      & \textbf{368 (133)} & \textbf{1501} & \textbf{834} & \textbf{2336} & \textbf{-8\%} 
      & \textbf{323 (69)} & \textbf{1258} & \textbf{762} & \textbf{2020} & \textbf{8\%}\\
      \hline
      \hline
      \textbf{VCS\_550\_1} & 254 & 1723 & 842 & 2565
                           & 275 (9) & 1125 & 913 & 2038 & \textbf{21\%}
                           & 368 (137) & 2101 & 1084 & 3185 & \textbf{-24\%}
                           & 309 (61) & 1421 & 899 & 2320 & \textbf{10\%}\\
      \textbf{VCS\_550\_2} & 280 & 1347 & 894 & 2268
                           & 280 (18) & 801 & 971 & 1772 & \textbf{22\%}
                           & 393 (145) & 1583 & 1116 & 2699 & \textbf{-19\%}
                           & 324 (87) & 1289 & 976 & 2265 & \textbf{0\%}\\
      \textbf{VCS\_550\_3} & 286 & 1245 & 822 & 2067
                           & 269 (21) & 659 & 829 & 1488 & \textbf{28\%}
                           & 326 (103) & 960 & 943 & 1903 & \textbf{8\%}
                           & 302 (70) & 888 & 907 & 1795 & \textbf{13\%}\\
      \textbf{VCS\_550\_4} & 259 & 1870 & 863 & 2733
                           & 262 (13) & 980 & 895 & 1875 & \textbf{31\%}
                           & 406 (132) & 2130 & 1186 & 3316 & \textbf{-21\%}
                           & 294 (48) & 1207 & 954 & 2161 & \textbf{21\%}\\
      \textbf{VCS\_550\_5} & 199 & 782 & 774 & 1556
                           & 206 (14) & 472 & 788 & 1260 & \textbf{19\%}
                           & 242 (67) & 632 & 905 & 1537 & \textbf{1\%}
                           & 228 (54) & 608 & 858 & 1466 & \textbf{6\%}\\
      \hline
      \textbf{Average} 
      & \textbf{256} & \textbf{1399} & \textbf{839} & \textbf{2238} 
      & \textbf{258 (15)} & \textbf{807} & \textbf{879} & \textbf{1687} & \textbf{24\%} 
      & \textbf{347 (116)} & \textbf{1481} & \textbf{1047} & \textbf{2528} & \textbf{-11\%} 
      & \textbf{291 (64)} & \textbf{1083} & \textbf{919} & \textbf{2001} & \textbf{10\%}\\
      \hline
      \hline
      
      \textbf{VCS\_600\_1} & 290 & 2840 & 1229 & 4069
                           & 298 (32) & 1694 & 1235 & 2929 & \textbf{28\%}
                           & 442 (171) & 3416 & 1595 & 5011 & \textbf{-23\%}
                           & 320 (69) & 1972 & 1300 & 3272 & \textbf{20\%}\\
      \textbf{VCS\_600\_2} & 246 & 1286 & 948 & 2234
                           & 249 (20) & 815 & 972 & 1787 & \textbf{20\%}
                           & 319 (114) & 1327 & 1136 & 2463 & \textbf{-10\%}
                           & 331 (81) & 1324 & 1185 & 2509 & \textbf{-12\%}\\
      \textbf{VCS\_600\_3} & 279 & 2120 & 1398 & 3518
                           & 290 (11) & 1236 & 1344 & 2580 & \textbf{27\%}
                           & 392 (114) & 2208 & 1693 & 3901 & \textbf{-11\%}
                           & 373 (57) & 1937 & 1537 & 3474 & \textbf{1\%}\\
      \textbf{VCS\_600\_4} & 258 & 2103 & 1327 & 3430
                           & 272 (15) & 1286 & 1322 & 2608 & \textbf{24\%}
                           & 379 (104) & 2483 & 1618 & 4101 & \textbf{-20\%}
                           & 321 (48) & 1767 & 1475 & 3242 & \textbf{5\%}\\
      \textbf{VCS\_600\_5} & 247 & 1722 & 1186 & 2908
                           & 245 (16) & 905 & 1190 & 2095 & \textbf{28\%}
                           & 313 (96) & 1554 & 1395 & 2949 & \textbf{-1\%}
                           & 290 (54) & 1510 & 1315 & 2825 & \textbf{3\%}\\
      \hline
      \textbf{Average} 
      & \textbf{264} & \textbf{2014} & \textbf{1218} & \textbf{3232} 
      & \textbf{270 (18)} & \textbf{1187} & \textbf{1213} & \textbf{2400} & \textbf{25\%} 
      & \textbf{369 (119)} & \textbf{2198} & \textbf{1487} & \textbf{3685} & \textbf{-13\%} 
      & \textbf{327 (61)} & \textbf{1702} & \textbf{1362} & \textbf{3064} & \textbf{3\%}\\
      \hline
    \end{tabular}}
    \caption{VCSP results.}
    \label{tab:table-results-vcsp}
  \end{center}
\end{table}

According to the results, if we start by comparing Baseline-CG and ML-S, we can observe a reduction in computing time ranging from 23 to 27\%, mostly on the PP side. One can also notice that the highest reductions are obtained when the number of iterations is lower than the one with Baseline-CG. On average, the number of iterations is very close or slightly higher. Therefore, the ML selection allows a reduction of the size of the PP network without increasing considerably the number of iterations. We can also notice that the performance of the model is maintained across all groups of instances. Considering that the model has been trained only on 400-trip instances, this also shows the ability of the model to generalize to different instance sizes. 

For the random selection case (i.e., Random-S), from the negative gains obtained we can conclude that this selection is not very promising. This is reflected by the number of iterations performed that is higher than the one obtained using Baseline-CG, and also by the number of times the full network was used if we compare it to ML-S. On the other hand, the results obtained by Cost-S show that a selection based on the costs is more effective. One can consider that Cost-S is a selection based on a single feature, which is the cost, while ML-S uses several features at the same time to decide whether an arc is to be retained or not. The results lie between those of the two algorithms Random-S and ML-S, i.e., it is better than a completely random selection but worse than the one with a learned model. This is reflected by the total number of iterations, the number of times the full network $G$ was used and the reductions obtained.

\subsubsection{\label{subsubsec:case1_additional_exp}Additional experiments}

This section is dedicated to some additional experiments that have been conducted. Even though they did not give significantly better results, we think that they deserve to be highlighted here.

\begin{itemize}
	\item[$-$] The first idea has been previously mentioned, which consists of selecting only the arcs that are part of the columns in the optimal basis of the solved RMPs. In the training phase, this strategy gave slightly worse results, i.e., 75\% accuracy, and when integrated in the CG algorithm, the extra gain obtained by reducing further the network size is nullified by doing more iterations.
	
	\item[$-$] Another idea is to assign weights to the arcs according to the number of times they have been used in the generated columns. An accuracy of 78\% is obtained in the training phase, distributed as 86\% recall and 70\% TNR, but again, the results were not any better when integrated in CG.
	
	\item[$-$] We made a small modification to the initial algorithm which consists of skipping the label assignment at the first $n$ iterations, where $n$ is a given parameter. The idea is that some arcs can be interesting at the beginning of the optimization but never used at later iterations. Like the other changes, the results obtained were not really compelling.
	
	\item[$-$] In the same context of combining ML and CG, in our previous paper \citet{mouad-columnselection}, we were interested in reducing the computing time of the RMPs by selecting the most promising columns at each iteration. Reductions in computation time of up to 30\% were achievable. However, the developed approach was more intended for problems that take the majority of the time in solving the RMPs and not the PPs, e.g., the tests were performed on instances that take on average 75\% or more of the computing time in solving the RMPs. Nevertheless, it was interesting to check if column selection can yield an additional gain when used with ML-S.
	The results comparing ML-S and ML-S with column selection, i.e., ML-COL-S, are reported in Table \ref{tab:table-results-vcsp-selcol}.
	\\ \\
	From these results, we observe a small additional reduction in the average computing time per instance group ranging from 5\% to 14\%. 
\end{itemize}

\begin{table}
  \begin{center}
    \resizebox{.7\textwidth}{!}{
    \begin{tabular}{|r||c|r|r|r||c|r|r|r|r|}
      \hline
      \multirow{3}{*}{\textbf{Instance}} &   \multicolumn{4}{|c||}{\textbf{ML-S}} &  \multicolumn{4}{|c|}{\textbf{ML-COL-S}} &
      \multirow{3}{*}{\textbf{Gain}}\\
      \cline{2-9}
      & \multirow{2}{*}{\#Itr} & \multicolumn{3}{|c||}{Time (s)}
      & \multirow{2}{*}{\#Itr} & \multicolumn{3}{|c|}{Time (s)} &  
      \\
      \cline{3-5} \cline{7-9}
      & & PP & RMP & Total
      & & PP & RMP & Total & \\
      
      \hline\hline
      \textbf{VCS\_300\_1} & 146 (12) & 75 & 76 & 151
                           & 123 (10) & 64 & 70 & 134 & \textbf{11\%}\\
      \textbf{VCS\_300\_2} & 158 (8) & 106 & 91 & 197
                           & 140 (7) & 103 & 78 & 181 & \textbf{8\%}\\ 
      \textbf{VCS\_300\_3} & 154 (14) & 115 & 76 & 191
                           & 135 (10) & 101 & 66 & 167 & \textbf{13\%}\\
      \textbf{VCS\_300\_4} & 185 (12) & 166 & 81 & 247
                           & 161 (6) & 128 & 67 & 195 & \textbf{21\%}\\
      \textbf{VCS\_300\_5} & 191 (27) & 133 & 87 & 220
                           & 142 (9) & 102 & 76 & 178 & \textbf{19\%}\\                    
      \hline
      \textbf{Average} 
      & \textbf{167 (14)} & \textbf{119} & \textbf{82} & \textbf{201} 
      & \textbf{140 (8)} & \textbf{99} & \textbf{71} & \textbf{171} & \textbf{14\%}\\
      \hline
      \hline
      \textbf{VCS\_350\_1} & 143 (15) & 112 & 129 & 241
                           & 126 (7) & 106 & 125 & 231 & \textbf{4\%}\\
      \textbf{VCS\_350\_2} & 135 (9) & 85 & 119 & 204
                           & 124 (9) & 81 & 105 & 186 & \textbf{9\%}\\ 
      \textbf{VCS\_350\_3} & 170 (10) & 141 & 167 & 308
                           & 163 (7) & 139 & 156 & 295 & \textbf{4\%}\\
      \textbf{VCS\_350\_4} & 188 (7) & 173 & 174 & 347
                           & 179 (9) & 198 & 148 & 346 & \textbf{0\%}\\
      \textbf{VCS\_350\_5} & 234 (12) & 272 & 183 & 455
                           & 187 (9) & 234 & 155 & 389 & \textbf{15\%}\\                    
      \hline
      \textbf{Average} 
      & \textbf{174 (10)} & \textbf{157} & \textbf{154} & \textbf{311} 
      & \textbf{155 (8)} & \textbf{151} & \textbf{137} & \textbf{289} & \textbf{6\%}\\
      \hline
      \hline
      \textbf{VCS\_400\_1} & 177 (10) & 238 & 258 & 496
                           & 161 (7) & 229 & 232 & 461 & \textbf{7\%}\\
      \textbf{VCS\_400\_2} & 189 (8) & 255 & 258 & 513
                           & 175 (11) & 255 & 228 & 483 & \textbf{6\%}\\ 
      \textbf{VCS\_400\_3} & 163 (9) & 177 & 233 & 410
                           & 148 (7) & 176 & 203 & 379 & \textbf{8\%}\\
      \textbf{VCS\_400\_4} & 147 (12) & 126 & 190 & 316
                           & 122 (6) & 108 & 166 & 274 & \textbf{13\%}\\
      \textbf{VCS\_400\_5} & 180 (12) & 283 & 230 & 513
                           & 200 (7) & 316 & 236 & 552 & \textbf{-8\%}\\                    
      \hline
      \textbf{Average} 
      & \textbf{171 (10)} & \textbf{216} & \textbf{234} & \textbf{450} 
      & \textbf{161 (7)} & \textbf{216} & \textbf{213} & \textbf{429} & \textbf{5\%}\\
      \hline
      \hline
      \textbf{VCS\_450\_1} & 368 (14) & 809 & 636 & 1445
                           & 319 (10) & 723 & 513 & 1236 & \textbf{14\%}\\
      \textbf{VCS\_450\_2} & 241 (11) & 797 & 486 & 1283
                           & 232 (6) & 833 & 430 & 1263 & \textbf{2\%}\\ 
      \textbf{VCS\_450\_3} & 194 (10) & 319 & 391 & 710
                           & 173 (9) & 314 & 355 & 669 & \textbf{6\%}\\
      \textbf{VCS\_450\_4} & 234 (8) & 312 & 440 & 752
                           & 190 (9) & 280 & 349 & 629 & \textbf{16\%}\\
      \textbf{VCS\_450\_5} & 289 (21) & 976 & 564 & 1540
                           & 263 (8) & 918 & 489 & 1407 & \textbf{9\%}\\                    
      \hline
      \textbf{Average} 
      & \textbf{265 (12)} & \textbf{643} & \textbf{503} & \textbf{1146} 
      & \textbf{235 (8)} & \textbf{613} & \textbf{427} & \textbf{1040} & \textbf{9\%}\\
      \hline
      \hline
      \textbf{VCS\_500\_1} & 304 (8) & 1030 & 902 & 1932
                           & 292 (7) & 986 & 888 & 1874 & \textbf{3\%}\\
      \textbf{VCS\_500\_2} & 270 (13) & 846 & 576 & 1422
                           & 259 (8) & 852 & 570 & 1422 & \textbf{0\%}\\ 
      \textbf{VCS\_500\_3} & 221 (10) & 453 & 611 & 1064
                           & 177 (7) & 416 & 550 & 966 & \textbf{9\%}\\
      \textbf{VCS\_500\_4} & 281 (8) & 1235 & 749 & 1984
                           & 247 (6) & 1134 & 670 & 1804 & \textbf{9\%}\\
      \textbf{VCS\_500\_5} & 340 (23) & 1024 & 656 & 1680
                           & 299 (11) & 916 & 604 & 1520 & \textbf{10\%}\\                    
      \hline
      \textbf{Average} 
      & \textbf{283 (12)} & \textbf{918} & \textbf{699} & \textbf{1616} 
      & \textbf{254 (7)} & \textbf{860} & \textbf{656} & \textbf{1517} & \textbf{6\%}\\
      \hline
      
    \end{tabular}}
    \caption{VCSP results with column selection.}
    \label{tab:table-results-vcsp-selcol}
  \end{center}
\end{table}

\section{\label{sec:case2}Application II: Vehicle routing problem with time windows}

The vehicle routing problem (VRP) is a well-known and classical problem that has been extensively studied in the literature. Given a fleet of vehicles, the VRP consists of constructing routes to serve geographically dispersed customers, while minimizing the travel costs and respecting the capacity of the vehicles. Each route must start at the depot, visit a set of clients and then return to the depot at the end of the tour. The variant with time windows (VRPTW) adds an additional complexity by introducing a restriction on the times during which the clients can be visited. Each client must be served during the associated time window by exactly one vehicle, and in case the vehicle arrives earlier, it has the possibility to wait until the client is available.

The most efficient and recent methods to solve this problem are based on CG, embedded in a branch-and-bound framework, with the addition of cuts to improve the bounds quality, giving rise to what is called ``\textit{branch-price-and-cut}" algorithms.
CG is used to solve the linear relaxation at each node of the tree, where the master problem is most likely a set partitioning problem and the PP is an ESPPRC solved using dynamic programming algorithms. The elementarity requirement of the PP makes the problem much more difficult to solve than the standard SPPRC, especially for large instances where several clients can be visited by one route. To overcome these difficulties, researchers have investigated several relaxations (see \citet{baldacci-ng-routes-2011}, \citet{irnich-spprc-k-cyc} and \citet{desaulniers-lessard-tabu-kpath}), seeking a good trade-off between the quality of the bound and the difficulty to solve the problem. Other studies have explored heuristic solutions (\citet{desaulniers-lessard-tabu-kpath}, \citet{fukasawa2006}), such as dominating on a subset of resources, limiting the number of labels kept at each node, or reducing the number of arcs in the network based on their reduced cost. Note that many of these techniques can be combined and used simultaneously, allowing a large reduction in computing time. For more details on the exact methods for solving the VRP, the survey of \citet{costa-contardo-desaulniers-2019} covers a wide range of methods and explores different variants of the problem. Similar to the recent formulations proposed in the literature, the model used in this work is a set partitioning problem and the PP is a SPPRC with 2-cycle elimination \citep{desrochers-vrptw-1992} solved using a labeling algorithm. 

In this section, we will first describe the network structure of the problem, before giving more details about the the data collection and the instances used. Finally, we will discuss the results obtained.

\subsection{\label{subsec:case2_network_structure}Network structure}

We consider the basic network structure of the VRPTW with a single depot, where only one type of arcs exists, representing the possible vehicle movements between the depot and the clients and between the clients. Therefore, all the arcs will be targeted by the selection strategy, with the exception of the ones entering and leaving the depot. For other problem variants and depending on the network structure, it would be possible to limit the selection to only a subset of arcs as we did for the VCSP.

\subsection{\label{subsec:case2_data_collection}Data collection, features and labels}

In this section, a few additional details specific to the VRPTW are given regarding the features described in Section \ref{sec:features}. Starting with the resource consumptions on the arcs, two resources are considered as features: time and capacity (i.e., the customer demand). For each arc $(i,j) \in A, i \neq s, j \neq t$, the lower and upper bounds of the time windows of both client nodes $i$ and $j$ are collected.
For the capacity resource, since the windows are the same for all clients (i.e., $[0, Q]$ where $Q$ is the vehicle capacity), their bounds are not included in the features.\\
Finally, the labels are assigned as described in Section \ref{sec:labels}, i.e., the label $1$ is assigned to the arcs that are part of the generated routes, and $0$ for the others. 

\subsection{\label{subsec:case2_instances}VRPTW instances}

The VRPTW instances used are based on the Gehring \& Homberger instances \citep{homberger-2005}. Three classes of instances are available: the $R$ instances, standing for ``Random", where the clients positions are randomly scattered, the $C$ instances, standing for ``Cluster", where clients are grouped together in clusters and thus have more chance to be served by the same vehicle, and finally the $RC$ instances that are a combination of both. Moreover, for each of the three classes, there are two categories numbered $1$ and $2$ ($R1$, $C1$, $RC1$ and $R2$, $C2$, $RC2$). The difference between the two is in the width of the time windows: the category $1$ instances tend to have larger time windows than those of category $2$, and are therefore more difficult to solve.

Since we are only interested in solving the linear relaxation of the problem, and that most of the $R1$, $C1$ and $RC1$ instances we considered, i.e., $200$ and $400$ clients, are actually very easy to solve and take only a few seconds, we decided to only use the $R2$, $C2$ and $RC2$ instances. Some of these instances have been slightly modified, more precisely, the width of the time windows has been slightly reduced so that the computing time becomes reasonable and not very high. This is because some instances take hours to be solved, especially when no heuristics are used during the data collection phase.

Regarding the size of the instances, we used those of $200$ and $400$ clients. Half of the 200-client instances were used for training and testing in the ML phase (e.g., from the three classes), while the other half and the 400-client instances were used for testing the model obtained when integrated in the CG algorithm.

\subsection{\label{subsec:case2_results}Computational results}

As for the VCSP, this section will be divided into three parts. The first part is devoted to the ML phase where additional details about the features and data collection are covered as well as the results obtained. The second part presents the results of the model integration in the CG process, as well as a comparison with other algorithms. In the third part, we include some additional experiments that are worth mentioning. All the experiments were performed on a Linux machine with an i7-8700 CPU @ 3.20GHz and 64GB of RAM.

\subsubsection{\label{subsubsec:case2_ml_results}Machine learning}

Because of the noticeable differences in the feature value range from one instance to another, especially between instances of different sizes, a preprocessing step is performed to scale and normalize the data of each instance individually. A random forest model is fit on the training data, and the hyperpameters are tuned through a cross validation approach. A neural network model was also tested, but without obtaining a significantly better accuracy, so we decided to retain the RF model as we did for the VCSP. The best hyperparameter values obtained are described in Table \ref{tab:vrptw-ml-hyperparameter}, whereas the results on the test set are reported in Table \ref{tab:vrptw-metrics}.

\begin{table}[t]
  \begin{center}
    
    \begin{tabular}{|l|c|} 
      \hline
      \textbf{Hyperparameter} & \textbf{Value} \\
      \hline
       Bootstrap & True\\
       Max depth & 5\\
       Max features & 5\\
       Min samples per leaf & 50\\
       Min samples per split & 100\\
       Number of trees & 500\\
       Class weights& Balanced\\
      \hline
    \end{tabular}
  \end{center}
  \caption{Hyperparameters values used in the training phase.}
  \label{tab:vrptw-ml-hyperparameter}
\end{table}

\begin{table}[t]
  \begin{center}
    
    \begin{tabular}{|l|c|} 
      \hline
      \textbf{Metric} & \textbf{Value} \\
      \hline
      
      Recall & 93\%\\
      $TNR$ & 87\%\\
      Balanced Accuracy & 90\%\\
      \hline
    \end{tabular}
  \end{center}
  \caption{Metrics of the model obtained for the VRPTW.}
  \label{tab:vrptw-metrics}
\end{table}

The hyperparameter values are slightly different from those used for the VCSP, more precisely the max\_depth and max\_features parameters. The class weights are used to balance the two classes as before, since there is only $14\%$ and $9\%$ of arcs with label 1 in the 200 and 400 instances, respectively. According to the results in Table \ref{tab:vrptw-metrics}, we notice a higher accuracy than the one obtained for the VCSP. The arcs to be selected are correctly predicted with $93\%$ accuracy (i.e., Recall), whereas the TNR value is 87\%, which means that the number of arcs selected by the model will not differ much from the number selected by the expert. 

Note that the VCSP and the VRPTW are two completely different problems and they are not supposed to give comparable accuracies, especially considering that the arcs to be selected play different roles in each problem.

\subsubsection{\label{subsubsec:case2_cg_results}CG with arc selection}

Now let us evaluate the ML model performance when integrated in the CG process. To do so, we have chosen to compare different algorithms described as follows:

\textbf{Baseline-CG}: This corresponds to the standard CG algorithm without arc selection where the complete network is used in all iterations.

\textbf{ML-S}: This is the CG algorithm with the arc selection activated as described in Algorithm \ref{alg:cg_ml_pricing}, but with a small adjustment. In fact, for the VRPTW we noticed that towards the end of the optimization, the reduced network fails to generate columns for several consecutive iterations. Thus, the algorithm ends up solving the PP with both networks at multiple iterations (i.e., with $G_r$ followed by $G$), which slows down the optimization. Instead, we decided to disable the use of the reduced network as soon as it fails to generate a column until the end of the optimization.

\textbf{RedCost-S}: This corresponds to the CG algorithm using a known heuristic pricing strategy that reduces the number of arcs in the network. It has been used by several researchers, and has shown to be very effective. For each node, excluding the depot nodes $s$ and $t$, the heuristic consists in setting a minimum number $N^{min}$ of incoming and outgoing arcs to keep. At each iteration, the arcs are sorted by their reduced cost, the $N^{min}$ incoming and outgoing arcs with the least cost are kept, while the others are removed from the network. It is also possible to define multiple values $N^{min}_1, N^{min}_2, \dots, \infty$ where $N^{min}_1 < N^{min}_2 < \cdots < \infty$. The PP tries first to generate columns with the parameter value $N^{min}_1$, if it fails it moves to the next parameter value, and so on. The values used in this experiment are $N^{min}_1=10$, $N^{min}_2=20$, $N^{min}_3=\infty$. Note that this algorithm performs the selection at each iteration while ML-S performs the selection only once before the beginning of the optimization.

The results obtained to compare the three algorithms are described in Table \ref{tab:table-results-vrptw-1}. For each algorithm, we report the same information as in Table \ref{tab:table-results-vcsp}. An additional column presents the time gain in percentage obtained by ML-S and RedCost-S in comparison to Baseline-CG. For ML-S, the number of iterations with the full network is indicated between parentheses.

\begin{table}
  \begin{center}
    \resizebox{\textwidth}{!}{
    \begin{tabular}{|c||c||c|c|c|c||c|c|c|c|c||c|c|c|c|c|}
      \hline
      \multirow{3}{*}{\textbf{Size}} & \multirow{3}{*}{\textbf{Instance}} &   \multicolumn{4}{|c||}{\textbf{Baseline-CG}} &  \multicolumn{5}{|c||}{\textbf{ML-S}} & 
      \multicolumn{5}{|c|}{\textbf{RedCost-S}} \\
      \cline{3-16}
      & & \multirow{2}{*}{\#Itr} & \multicolumn{3}{|c||}{Time (s)}
      & \multirow{2}{*}{\#Itr} & \multicolumn{3}{|c|}{Time (s)} & \multirow{2}{*}{\textbf{Gain}} 
      & \multirow{2}{*}{\#Itr} & \multicolumn{3}{|c|}{Time (s)} & \multirow{2}{*}{\textbf{Gain}} \\
      \cline{4-6} \cline{8-10} \cline{13-15} 
     & & & PP & RMP & Total
      & & PP & RMP & Total & 
      & & PP & RMP & Total & \\
      
      \hline \hline
     \multirow{18}{*}{\rotatebox[origin=c]{90}{\textbf{200 clients}}} & \textbf{R2\_2\_1} & 197 & 173 & 8 & 181
                           & 259 (165) & 52 & 8 & 54 & \textbf{68\%}
                           & 439 & 23 & 12 & 34 & \textbf{81\%}\\
      & \textbf{R2\_2\_2} & 159 & 312 & 5 & 317
                           & 211 (150) & 99 & 6 & 105 & \textbf{67\%}
                           & 179 & 25 & 5 & 29 & \textbf{91\%}\\
      & \textbf{R2\_2\_3} & 198 & 128 & 7 & 135
                           & 301 (233) & 36 & 8 & 44 & \textbf{68\%}
                           & 441 & 20 & 8 & 28 & \textbf{79\%}\\
      & \textbf{R2\_2\_4} & 114 & 210 & 2 & 212
                           & 152 (101) & 69 & 3 & 72 & \textbf{66\%}
                           & 152 & 34 & 2 & 37 & \textbf{83\%}\\
      & \textbf{R2\_2\_5} & 239 & 238 & 10 & 248
                           & 275 (166) & 44 & 9 & 53 & \textbf{79\%}
                           & 367 & 22 & 11 & 33 & \textbf{87\%}\\
      \cline{2-16}
      & \textbf{Average} 
      & \textbf{181} & \textbf{216} & \textbf{6} & \textbf{233} 
      & \textbf{240 (163)} & \textbf{59} & \textbf{7} & \textbf{65} & \textbf{70\%} 
      & \textbf{316} & \textbf{25} & \textbf{7} & \textbf{32} & \textbf{84\%}\\
      \cline{2-16}
      \noalign{\vskip-2\tabcolsep \vskip-3\arrayrulewidth \vskip\doublerulesep} \\
      \cline{2-16}
      
      & \textbf{RC2\_2\_1} & 204 & 203 & 5 & 208
                           & 243 (172) & 123 & 5 & 128 & \textbf{38\%}
                           & 374 & 32 & 6 & 38 & \textbf{82\%}\\
      & \textbf{RC2\_2\_2} & 119 & 138 & 3 & 141
                           & 156 (107) & 49 & 3 & 52 & \textbf{63\%}
                           & 218 & 17 & 3 & 20 & \textbf{86\%}\\
      & \textbf{RC2\_2\_3} & 243 & 211 & 7 & 218
                           & 314 (214) & 82 & 9 & 91 & \textbf{59\%}
                           & 852 & 40 & 11 & 51 & \textbf{77\%}\\
      & \textbf{RC2\_2\_4} & 198 & 189 & 7 & 196
                           & 239 (166) & 150 & 7 & 157 & \textbf{20\%}
                           & 609 & 34 & 13 & 47 & \textbf{76\%}\\
      & \textbf{RC2\_2\_5} & 230 & 295 & 6 & 301
                           & 278 (192) & 112 & 6 & 118 & \textbf{61\%}
                           & 313 & 32 & 6 & 39 & \textbf{87\%}\\
      \cline{2-16}
      & \textbf{Average} 
      & \textbf{199} & \textbf{207} & \textbf{6} & \textbf{213} 
      & \textbf{246 (170)} & \textbf{103} & \textbf{6} & \textbf{109} & \textbf{48\%} 
      & \textbf{473} & \textbf{31} & \textbf{8} & \textbf{39} & \textbf{81\%}\\
      \cline{2-16}
      \noalign{\vskip-2\tabcolsep \vskip-3\arrayrulewidth \vskip\doublerulesep} \\
      \cline{2-16}
      & \textbf{C2\_2\_1} & 288 & 177 & 8 & 185
                           & 398 (268) & 143 & 10 & 153 & \textbf{18\%}
                           & 760 & 72 & 11 & 83 & \textbf{55\%}\\
      & \textbf{C2\_2\_2} & 149 & 353 & 5 & 358
                           & 195 (126) & 255 & 5 & 260 & \textbf{27\%}
                           & 316 & 121 & 6 & 126 & \textbf{65\%}\\
      & \textbf{C2\_2\_3} & 211 & 190 & 6 & 196
                           & 359 (176) & 146 & 11 & 157 & \textbf{20\%}
                           & 272 & 54 & 7 & 61 & \textbf{69\%}\\
      & \textbf{C2\_2\_4} & 164 & 267 & 6 & 273
                           & 351 (163) & 230 & 9 & 239 & \textbf{12\%}
                           & 181 & 65 & 5 & 70 & \textbf{74\%}\\
      & \textbf{C2\_2\_5} & 158 & 469 & 6 & 475
                           & 420 (161) & 440 & 11 & 451 & \textbf{5\%}
                           & 237 & 140 & 6 & 146 & \textbf{69\%}\\
      \cline{2-16}
      & \textbf{Average} 
      & \textbf{194} & \textbf{291} & \textbf{6} & \textbf{297} 
      & \textbf{345 (178)} & \textbf{243} & \textbf{9} & \textbf{252} & \textbf{16\%} 
      & \textbf{353} & \textbf{90} & \textbf{7} & \textbf{97} & \textbf{66\%}\\
      \cline{1-16}
      \noalign{\vskip-2\tabcolsep \vskip-3\arrayrulewidth \vskip\doublerulesep} \\
      \cline{1-16}
    \multirow{28}{*}{\rotatebox[origin=c]{90}{\textbf{400 clients}}} & \textbf{R2\_4\_1} & 397 & 643 & 143 & 786
                           & 447 (239) & 154 & 104 & 258 & \textbf{67\%}
                           & 373 & 66 & 119 & 185 & \textbf{76\%}\\
    & \textbf{R2\_4\_2} & 248 & 318 & 85 & 403
                           & 330 (165) & 90 & 76 & 166 & \textbf{59\%}
                           & 289 & 45 & 84 & 129 & \textbf{68\%}\\
    & \textbf{R2\_4\_3} & 202 & 1058 & 42 & 1100
                           & 252 (188) & 333 & 47 & 380 & \textbf{65\%}
                           & 235 & 97 & 37 & 134 & \textbf{88\%}\\
    & \textbf{R2\_4\_4} & 149 & 734 & 21 & 755
                           & 194 (112) & 203 & 23 & 226 & \textbf{70\%}
                           & 239 & 142 & 20 & 162 & \textbf{79\%}\\
    & \textbf{R2\_4\_5} & 355 & 316 & 131 & 447
                           & 439 (244) & 124 & 107 & 231 & \textbf{48\%}
                           & 361 & 51 & 127 & 178 & \textbf{60\%}\\
    & \textbf{R2\_4\_6} & 240 & 293 & 76 & 369
                           & 331 (162) & 75 & 67 & 142 & \textbf{62\%}
                           & 277 & 44 & 78 & 122 & \textbf{67\%}\\
    & \textbf{R2\_4\_7} & 245 & 629 & 47 & 676
                           & 303 (179) & 227 & 50 & 277 & \textbf{59\%}
                           & 218 & 81 & 44 & 125 & \textbf{82\%}\\
    & \textbf{R2\_4\_8} & 187 & 632 & 25 & 657
                           & 203 (116) & 243 & 27 & 270 & \textbf{59\%}
                           & 242 & 144 & 23 & 167 & \textbf{75\%}\\
    & \textbf{R2\_4\_9} & 334 & 870 & 138 & 1008
                           & 472 (253) & 277 & 111 & 388 & \textbf{62\%}
                           & 322 & 101 & 116 & 217 & \textbf{78\%}\\
    & \textbf{R2\_4\_10} & 308 & 988 & 113 & 1101
                           & 461 (257) & 521 & 102 & 623 & \textbf{43\%}
                           & 380 & 154 & 108 & 262 & \textbf{76\%}\\
    \cline{2-16}
    & \textbf{Average} 
      & \textbf{267} & \textbf{648} & \textbf{82} & \textbf{730} 
      & \textbf{343 (191)} & \textbf{255} & \textbf{71} & \textbf{296} & \textbf{59\%} 
      & \textbf{294} & \textbf{93} & \textbf{76} & \textbf{168} & \textbf{75\%}\\
    \cline{2-16}
      \noalign{\vskip-2\tabcolsep \vskip-3\arrayrulewidth \vskip\doublerulesep} \\
      \cline{2-16}
    & \textbf{RC2\_4\_1} & 392 & 281 & 123 & 404
                           & 555 (288) & 144 & 116 & 260 & \textbf{36\%}
                           & 390 & 45 & 102 & 147 & \textbf{64\%}\\
    & \textbf{RC2\_4\_2} & 301 & 282 & 65 & 347
                           & 434 (216) & 108 & 70 & 178 & \textbf{49\%}
                           & 295 & 37 & 67 & 104 & \textbf{70\%}\\
    & \textbf{RC2\_4\_3} & 202 & 1224 & 29 & 1253
                           & 293 (151) & 288 & 33 & 321 & \textbf{74\%}
                           & 301 & 108 & 27 & 135 & \textbf{89\%}\\
    & \textbf{RC2\_4\_4} & 155 & 361 & 17 & 378
                           & 231 (125) & 185 & 22 & 207 & \textbf{45\%}
                           & 183 & 81 & 16 & 97 & \textbf{74\%}\\
    & \textbf{RC2\_4\_5} & 797 & 337 & 198 & 535
                           & 1301 (848) & 271 & 199 & 470 & \textbf{12\%}
                           & 1164 & 84 & 186 & 270 & \textbf{50\%}\\
    & \textbf{RC2\_4\_6} & 370 & 622 & 98 & 720
                           & 543 (295) & 229 & 101 & 330 & \textbf{54\%}
                           & 434 & 98 & 118 & 216 & \textbf{70\%}\\
    & \textbf{RC2\_4\_7} & 340 & 615 & 87 & 702
                           & 510 (235) & 264 & 92 & 356 & \textbf{49\%}
                           & 415 & 100 & 98 & 198 & \textbf{72\%}\\
    & \textbf{RC2\_4\_8} & 373 & 629 & 74 & 703
                           & 478 (254) & 245 & 73 & 318 & \textbf{55\%}
                           & 435 & 90 & 72 & 162 & \textbf{77\%}\\
    & \textbf{RC2\_4\_9} & 305 & 460 & 83 & 543
                           & 417 (263) & 188 & 91 & 279 & \textbf{49\%}
                           & 371 & 83 & 91 & 174 & \textbf{68\%}\\
    & \textbf{RC2\_4\_10} & 411 & 219 & 101 & 320
                           & 458 (347) & 128 & 89 & 217 & \textbf{32\%}
                           & 438 & 40 & 91 & 131 & \textbf{59\%}\\
    
    \cline{2-16}
      & \textbf{Average} 
      & \textbf{365} & \textbf{503} & \textbf{88} & \textbf{591} 
      & \textbf{522 (302)} & \textbf{205} & \textbf{89} & \textbf{294} & \textbf{46\%} 
      & \textbf{443} & \textbf{77} & \textbf{87} & \textbf{163} & \textbf{69\%}\\
    \cline{2-16}
      \noalign{\vskip-2\tabcolsep \vskip-3\arrayrulewidth \vskip\doublerulesep} \\
      \cline{2-16}
    & \textbf{C2\_4\_1} & 599 & 373 & 136 & 509
                           & 1125 (660) & 424 & 178 & 602 & \textbf{-18\%}
                           & 966 & 144 & 183 & 327 & \textbf{36\%}\\
    & \textbf{C2\_4\_2} & 270 & 214 & 74 & 288
                           & 484 (237) & 197 & 78 & 275 & \textbf{5\%}
                           & 445 & 71 & 88 & 159 & \textbf{45\%}\\
    & \textbf{C2\_4\_3} & 254 & 543 & 50 & 593
                           & 342 (180) & 496 & 58 & 554 & \textbf{7\%}
                           & 484 & 176 & 55 & 231 & \textbf{61\%}\\
    & \textbf{C2\_4\_4} & 153 & 574 & 26 & 600
                           & 408 (121) & 588 & 52 & 640 & \textbf{-7\%}
                           & 332 & 267 & 30 & 297 & \textbf{51\%}\\
    & \textbf{C2\_4\_5} & 261 & 632 & 62 & 694
                           & 714 (222) & 637 & 93 & 730 & \textbf{-5\%}
                           & 417 & 194 & 63 & 257 & \textbf{63\%}\\
    \cline{2-16}
    & \textbf{Average} 
      & \textbf{307} & \textbf{467} & \textbf{69} & \textbf{536} 
      & \textbf{614 (284)} & \textbf{468} & \textbf{91} & \textbf{560} & \textbf{-4\%} 
      & \textbf{528} & \textbf{170} & \textbf{83} & \textbf{254} & \textbf{51\%}\\
      \hline
    \end{tabular}}
    \caption{VRPTW results.}
    \label{tab:table-results-vrptw-1}
  \end{center}
\end{table}

\begin{figure}
     \centering
     \begin{subfigure}[b]{0.31\textwidth}
         \centering
         \includegraphics[width=\textwidth]{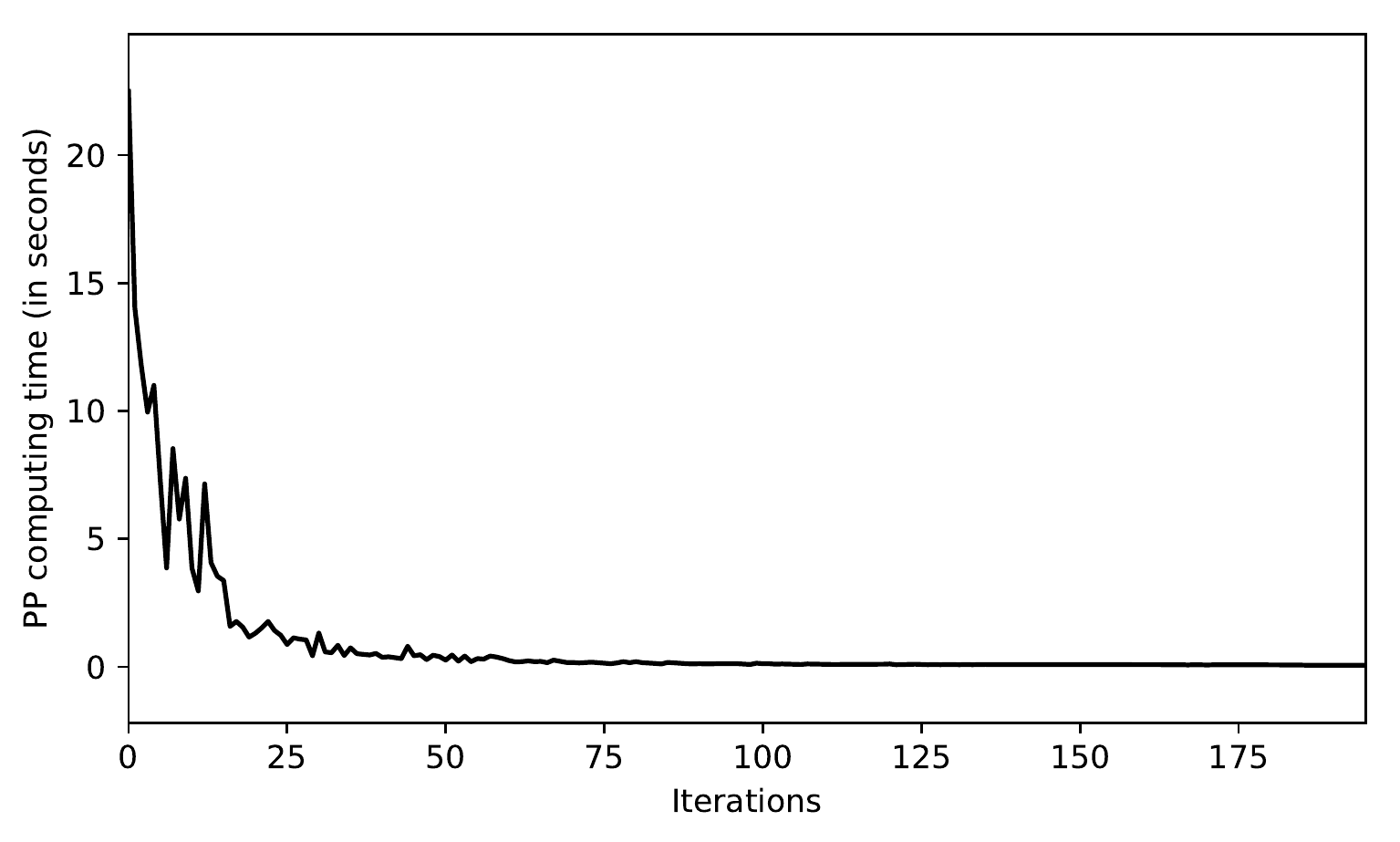}
         \caption{Computing time per iteration}
         \label{fig:baseline_computing}
     \end{subfigure}
     \hfill
     \begin{subfigure}[b]{0.31\textwidth}
         \centering
         \includegraphics[width=\textwidth]{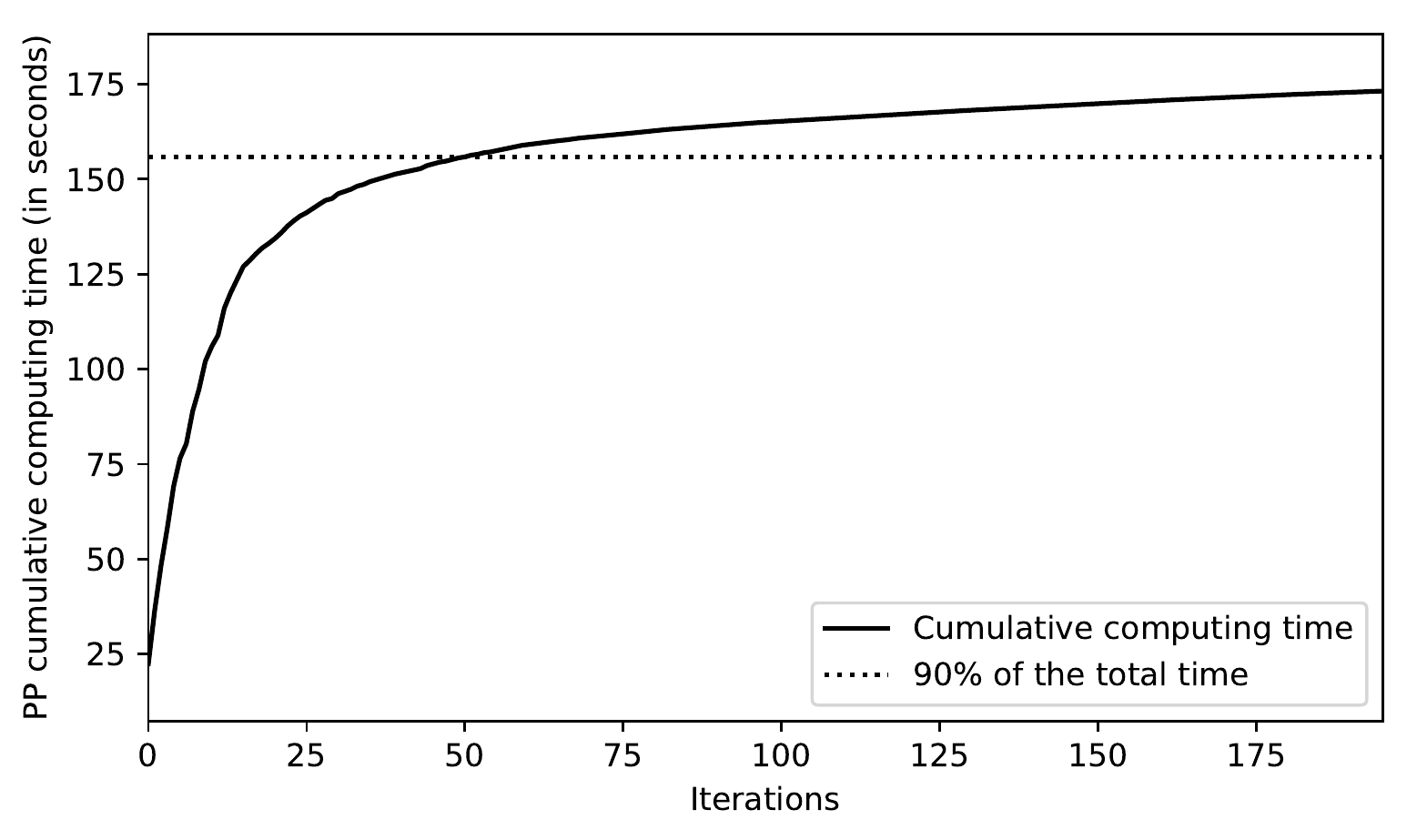}
         \caption{Cumulative computing time}
         \label{fig:baseline_computing_cumulative}
     \end{subfigure}
     \hfill
     \begin{subfigure}[b]{0.31\textwidth}
         \centering
         \includegraphics[width=\textwidth]{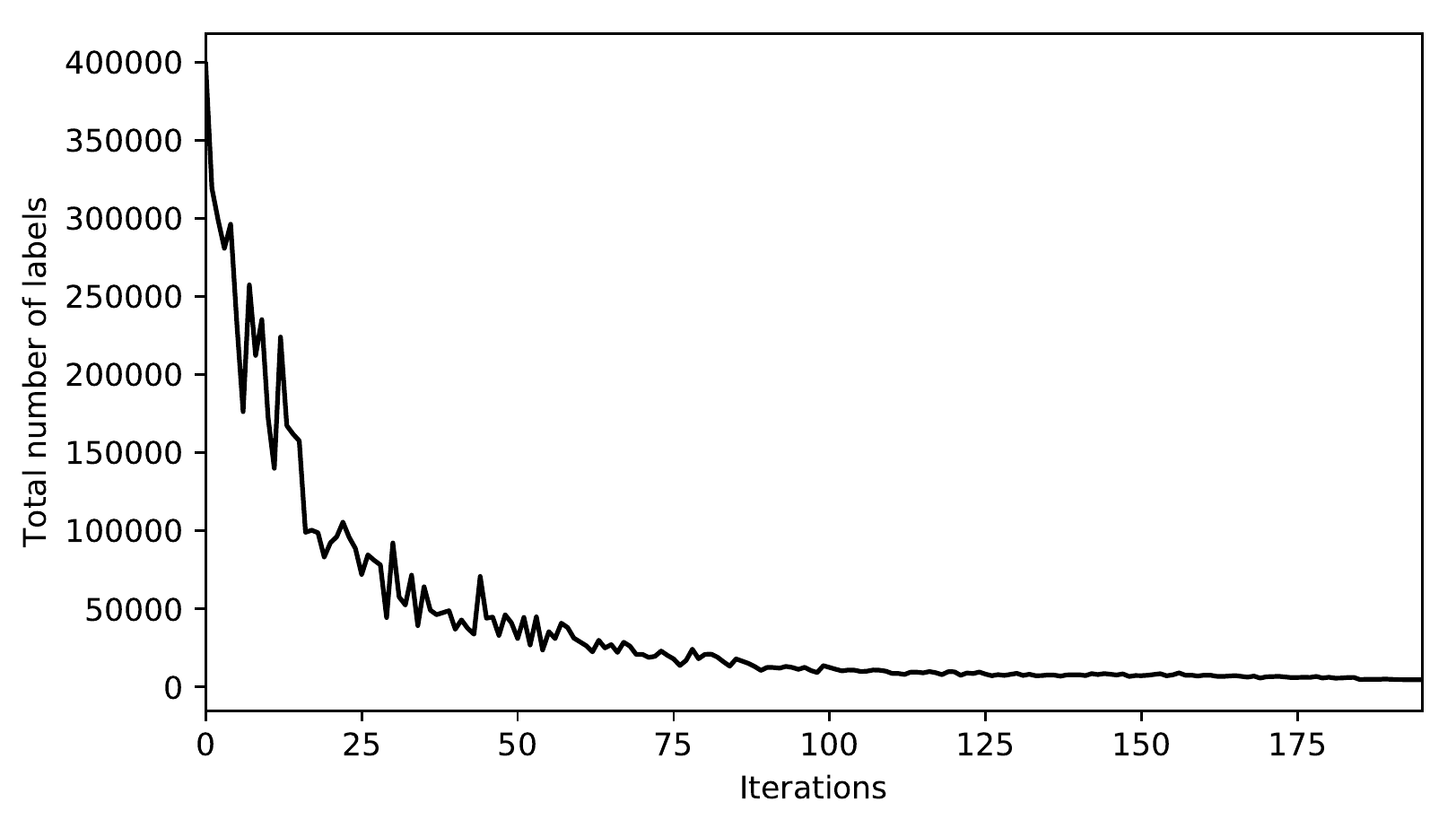}
         \caption{Number of labels per iteration}
         \label{fig:baseline_labels}
     \end{subfigure}
        \caption{PP computing time and number of labels with Baseline-CG.}
        \label{fig:baseline_results}
\end{figure}

Starting with Baseline-CG, we can notice that most of the computing time is spent on the PP side. If we take a closer look at the time per iteration of one of the instances (R2\_2\_1 as an example), as shown in Figure \ref{fig:baseline_computing}, the PP computing time decreases very rapidly, starting at $22$ seconds at the first iteration and reaching less than $0.1$ second after a hundred iterations.  In fact, according to the cumulative time shown in Figure \ref{fig:baseline_computing_cumulative}, we can see that $90\%$ of the computing time is spent in the first $50$ iterations, which is just a quarter of the total number of iterations. The progressive decrease of computing time can be explained by looking at the number of labels created during the resolution of the PP as shown in Figure \ref{fig:baseline_labels}. This shows that more labels are dominated as we progress in the optimization, and with less labels there is less processing and therefore less computing time. 

Moving on to the ML-S results, we can notice a higher number of iterations compared to Baseline-CG. The number of iterations with the full network (i.e., in parentheses) is also quite high, representing more than $50\%$ of the total number. Nevertheless, a significant reduction of computing time on the PP side is obtained. Since most of the iterations using the full network are performed in the second half of the optimization (i.e., when the use of the reduced network is disabled), they are not time consuming as previously shown in Figure \ref{fig:baseline_results}. We can see a significant gain on the $R2$ instances. However, this gain decreases on the instances where the clients are more clustered (i.e., RC2 and C2) and also when moving from the $200$ instances to the $400$ ones. If we compare the increase in the number of iterations between ML-S and Baseline-CG for the 200-client instances, we can notice an increase of $32\%$ for $R2$, $24\%$ for $RC2$ and $73\%$ for $C2$. It is true that the $C2$ instances seem to perform way more iterations than the other groups, and this may be one reason explaining the drop in performance. However, we see that the $RC2$ instances have a lower increase in the number of iterations than $R2$ and yet their results are less good.

\begin{figure}
     \centering
     \begin{subfigure}[b]{0.31\textwidth}
         \centering
         \includegraphics[width=\textwidth]{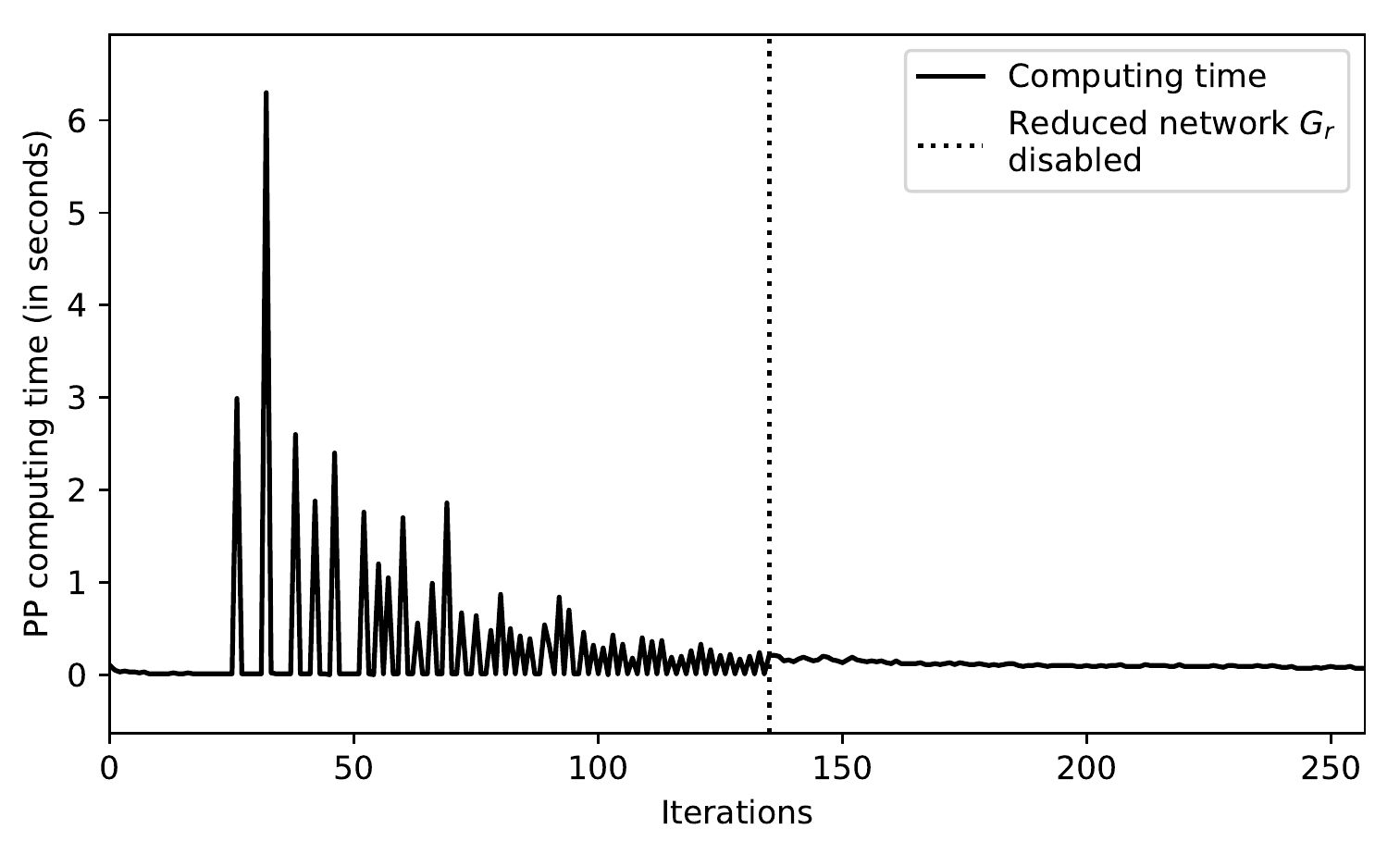}
         \caption{Computing time per iteration}
         \label{fig:r2_computing}
     \end{subfigure}
     \hfill
     \begin{subfigure}[b]{0.31\textwidth}
         \centering
         \includegraphics[width=\textwidth]{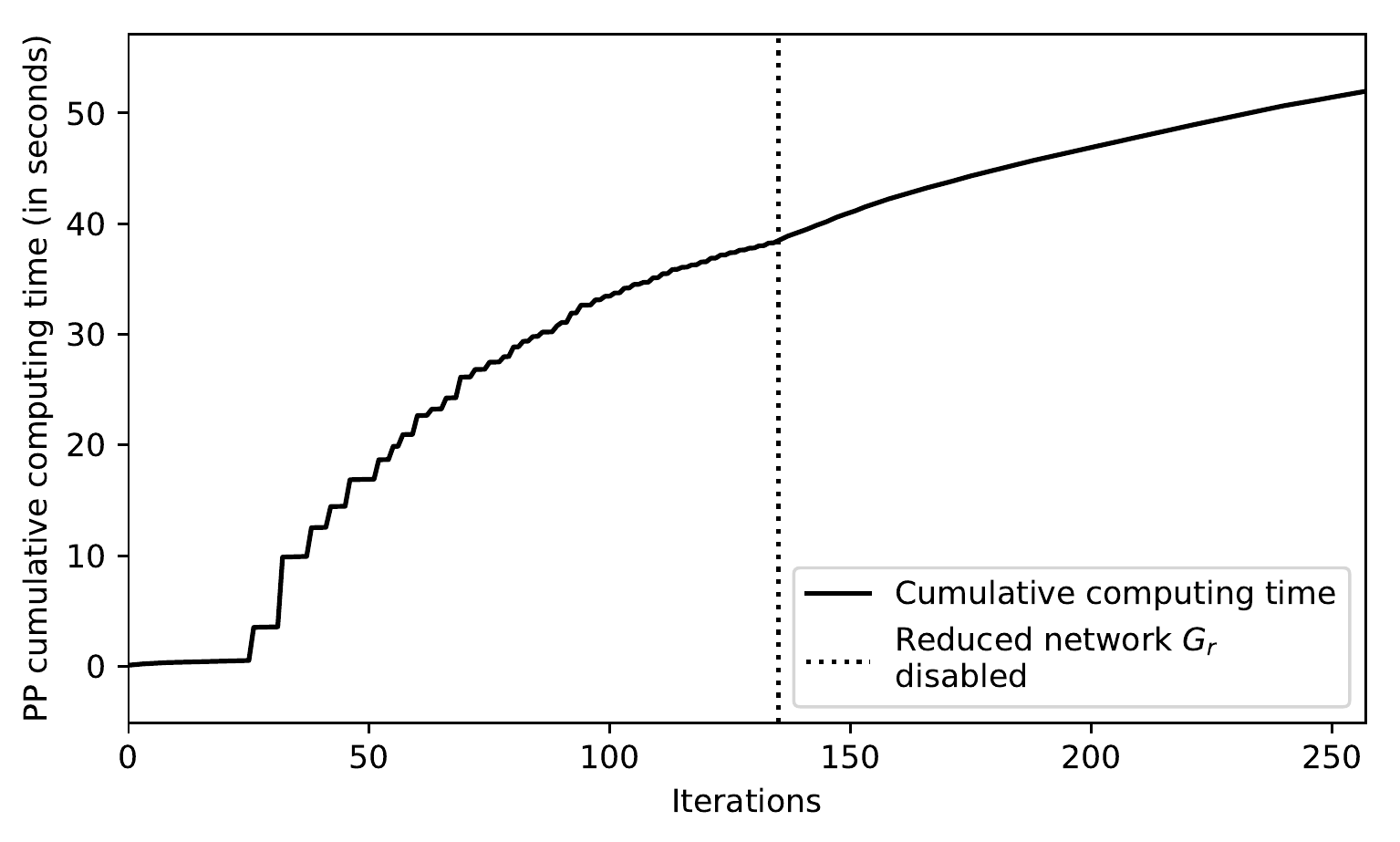}
         \caption{Cumulative computing time}
         \label{fig:r2_computing_cumulative}
     \end{subfigure}
     \hfill
     \begin{subfigure}[b]{0.31\textwidth}
         \centering
         \includegraphics[width=\textwidth]{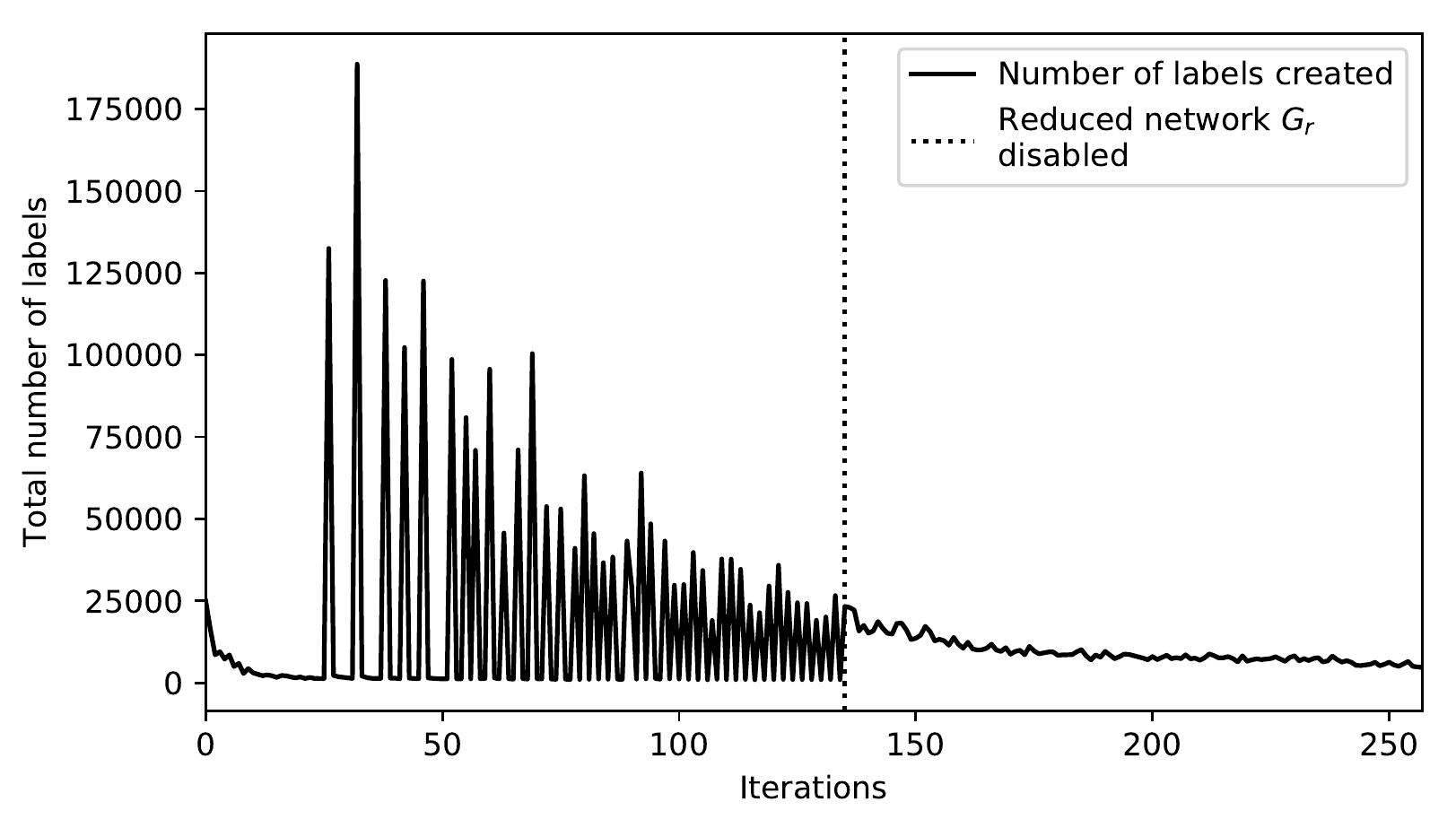}
         \caption{Number of labels per iteration}
         \label{fig:r2_labels}
     \end{subfigure}
        \caption{The PP computing time and number of labels with ML-S for instance R2\_2\_1.}
        \label{fig:r2_results}
    \vspace*{7mm}
     \centering
     \begin{subfigure}[b]{0.31\textwidth}
         \centering
         \includegraphics[width=\textwidth]{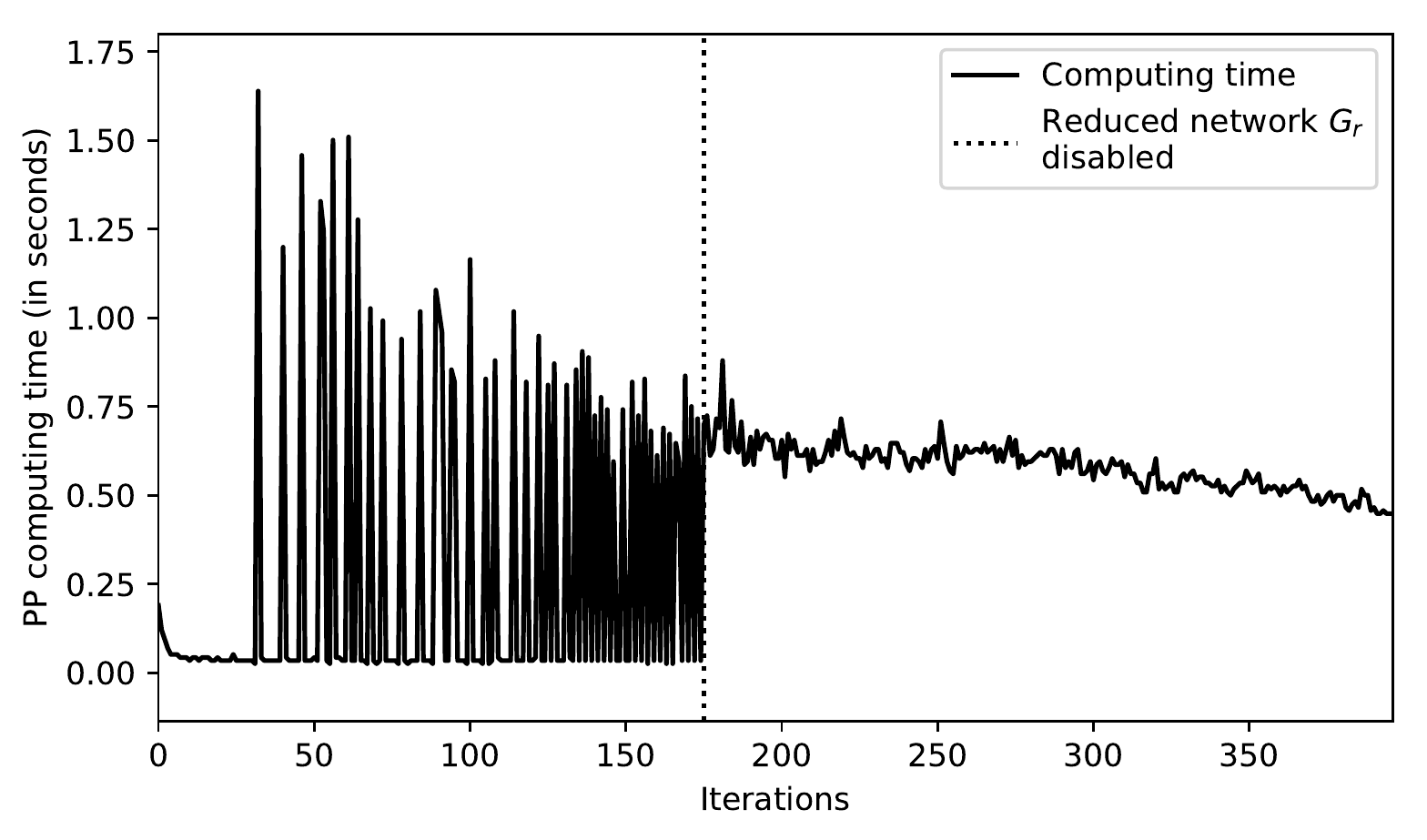}
         \caption{Computing time per iteration}
         \label{fig:c2_computing}
     \end{subfigure}
     \hfill
     \begin{subfigure}[b]{0.31\textwidth}
         \centering
         \includegraphics[width=\textwidth]{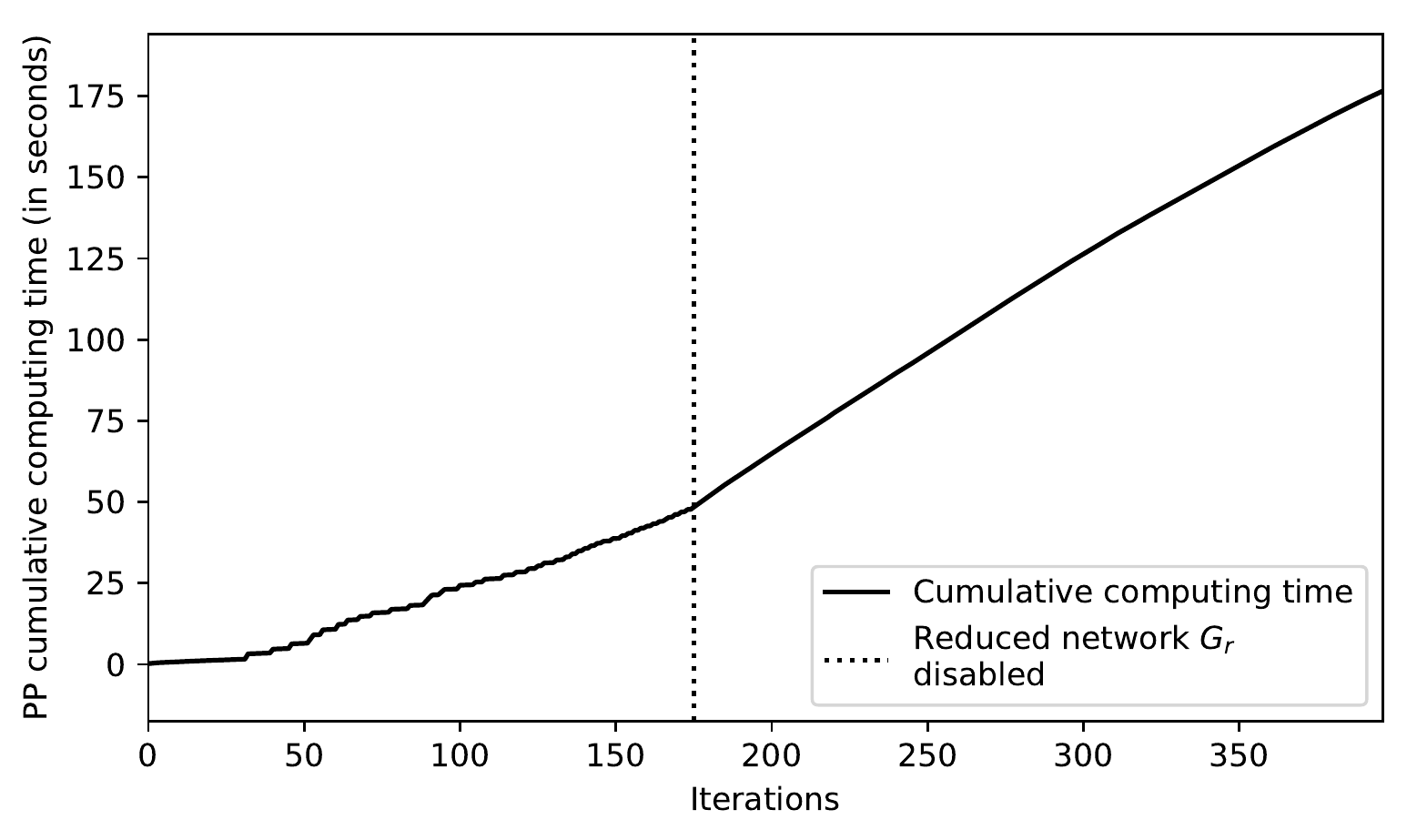}
         \caption{Cumulative computing time}
         \label{fig:c2_computing_cumulative}
     \end{subfigure}
     \hfill
     \begin{subfigure}[b]{0.31\textwidth}
         \centering
         \includegraphics[width=\textwidth]{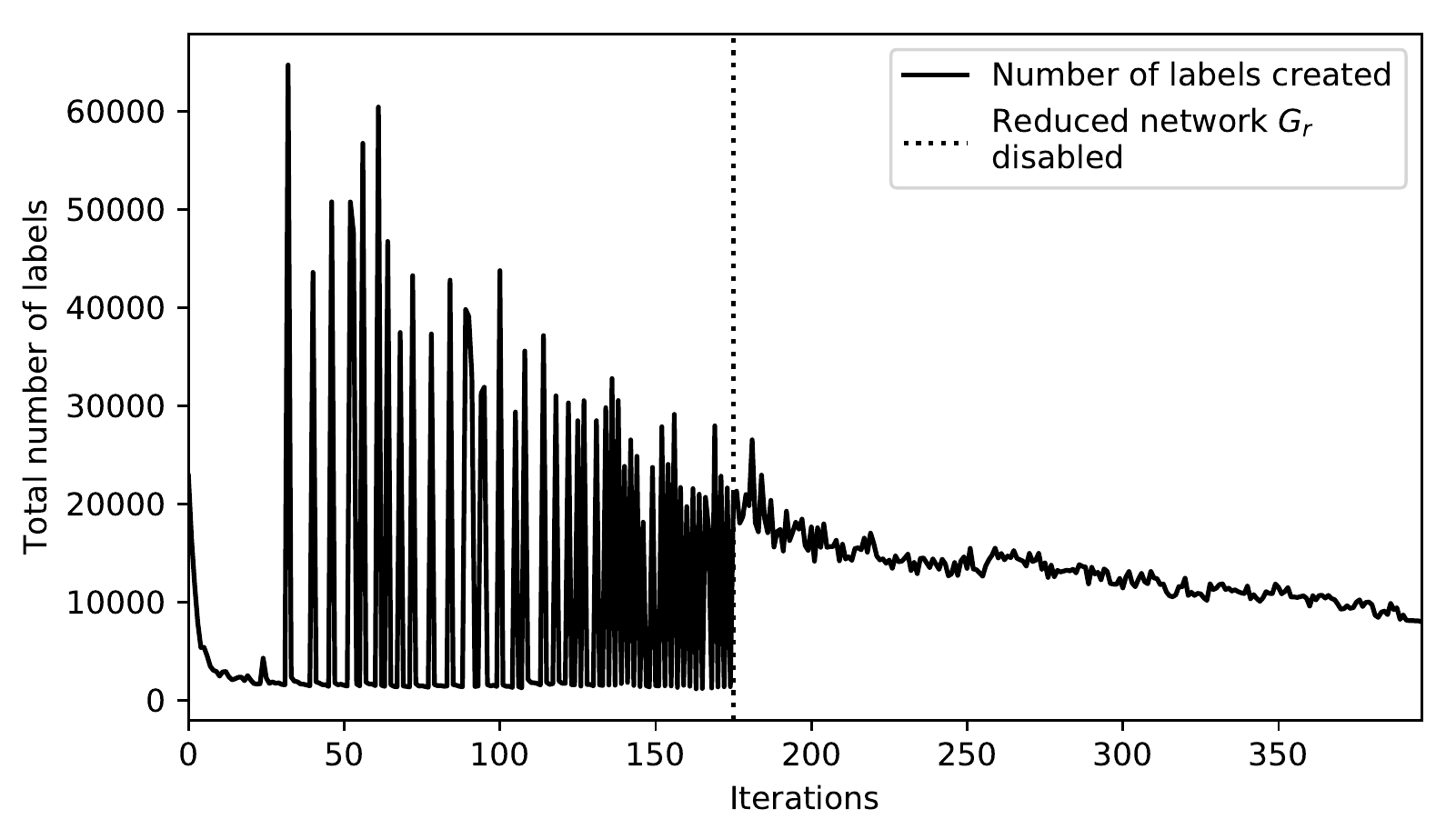}
         \caption{Number of labels per iteration}
         \label{fig:c2_labels}
     \end{subfigure}
        \caption{The PP computing time and number of labels with ML-S for instance C2\_2\_1.}
        \label{fig:c2_results}
\end{figure}

Let us consider that the optimization goes through two stages when using ML-S. The first stage starts from the beginning until the reduced network fails to generate any columns (i.e., the reduced network is thus disabled), and the second stage comes afterwards until the end of the optimization. By taking an $R2$ instance and comparing the computing time per iteration during the two stages, Figure \ref{fig:r2_computing} shows that there are multiple peaks of computing time during the first stage, which occur when we switch to the full network (i.e., when the reduced network does not generate enough columns). We can also see that the peaks heights decreases rapidly as we progress in the optimization. In the second stage (e.g., starting at the dashed line), the computing time is not too high, and it is also more stable since we are only using the full network. Figures \ref{fig:r2_computing_cumulative} and \ref{fig:r2_labels} give some additional information, namely the cumulative time and the number of labels created, respectively. According to the cumulative time, we can see a step line in the first stage where the computing time increases significantly each time we switch to the full network, whereas in the second stage the line is almost linear. Notice that the total time spent on the first stage is larger than the second one. On the other hand, the number of labels follows exactly the same trend depicted in Figure \ref{fig:r2_computing}, which makes sense.

If we try to visualize the same data for a $C2$ instance in Figure \ref{fig:c2_results} and compare it to the previously seen data for an $R2$ instance (i.e., Figure \ref{fig:r2_results}), some differences can be noticed. We can observe that the computing time per iteration during the first stage decreases more slowly and is still quite high when we reach the second stage. The same thing can be noticed for the number of labels created. Looking at the cumulative time, we can see that the second stage takes more time than the first and the line is also steeper.

Finally, for the RedCost-S algorithm, we can notice that the results are significantly better than ML-S. A drop in performance is observed when dealing with the more clustered instances, but it is not as aggressive as ML-S. The number of iterations for the 200-client instances is larger or very close to the one of ML-S, but for the 400-client instances it is less on average.

Now the question is whether we can achieve better results than those obtained with RedCost-S. To answer that, we tried two additional strategies:

\textbf{ML-RedCost-S}: This corresponds to the combination of the two strategies ML-S and RedCost-S. In other words, the arc reduction according to the reduced cost is performed at each iteration on the two graphs $G$ and $G_r$ depending on which one is used.

\textbf{ML-RedCost-S-2}: This corresponds to an ML model learned on the data generated using RedCost-S. This means that RedCost-S is considered the new expert. By going back to the first phase, i.e., data collection, the label 1 is assigned to the arcs that are part of the columns generated using RedCost-S, and 0 for the others. The training is performed on the new data in order to obtain a new model, which is then used to build the reduced network. 

For both strategies, the same values of $N^{min}_1=10$, $N^{min}_2=20$, $N^{min}_3=\infty$ are used. The results comparing the two strategies to RedCost-S are reported in Table \ref{tab:table-results-vrptw-2}, providing the same information as in Table \ref{tab:table-results-vrptw-1} except that the time gains are measured with respect to the time obtained by RedCost-S.

\begin{table}
  \begin{center}
    \resizebox{\textwidth}{!}{
    \begin{tabular}{|c||c||c|c|c|c||c|c|c|c|c||c|c|c|c|c|}
      \hline
      \multirow{3}{*}{\textbf{Size}} & \multirow{3}{*}{\textbf{Instance}} &   \multicolumn{4}{|c||}{\textbf{RedCost-S}} &  \multicolumn{5}{|c||}{\textbf{ML-RedCost-S}} & 
      \multicolumn{5}{|c|}{\textbf{ML-RedCost-S-2}} \\
      \cline{3-16}
      & & \multirow{2}{*}{\#Itr} & \multicolumn{3}{|c||}{Time (s)}
      & \multirow{2}{*}{\#Itr} & \multicolumn{3}{|c|}{Time (s)} & \multirow{2}{*}{\textbf{Gain}} 
      & \multirow{2}{*}{\#Itr} & \multicolumn{3}{|c|}{Time (s)} & \multirow{2}{*}{\textbf{Gain}} \\
      \cline{4-6} \cline{8-10} \cline{13-15} 
     & & & PP & RMP & Total
      & & PP & RMP & Total & 
      & & PP & RMP & Total & \\
      
      \hline \hline
     \multirow{18}{*}{\rotatebox[origin=c]{90}{\textbf{200 clients}}} & \textbf{R2\_2\_1} & 439 & 23 & 12 & 34
                           & 578 (403) & 13 & 12 & 24 & \textbf{31\%}
                           & 420 (235) & 8 & 10 & 18 & \textbf{48\%}\\
      & \textbf{R2\_2\_2} & 179 & 25 & 5 & 29
                           & 309 (127) & 13 & 7 & 21 & \textbf{30\%}
                           & 343 (158) & 16 & 8 & 24 & \textbf{17\%}\\
      & \textbf{R2\_2\_3} & 441 & 20 & 8 & 28
                           & 406 (255) & 8 & 9 & 17 & \textbf{38\%}
                           & 378 (219) & 8 & 8 & 16 & \textbf{44\%}\\
      & \textbf{R2\_2\_4} & 152 & 34 & 2 & 37
                           & 217 (98) & 8 & 3 & 11 & \textbf{69\%}
                           & 245 (113) & 9 & 4 & 12 & \textbf{66\%}\\
      & \textbf{R2\_2\_5} & 367 & 22 & 11 & 33
                           & 430 (244) & 11 & 12 & 23 & \textbf{32\%}
                           & 444 (245) & 11 & 12 & 23 & \textbf{32\%}\\
      \cline{2-16}
      & \textbf{Average} 
      & \textbf{316} & \textbf{25} & \textbf{7} & \textbf{32} 
      & \textbf{388 (225)} & \textbf{11} & \textbf{8} & \textbf{19} & \textbf{40\%} 
      & \textbf{366 (194)} & \textbf{10} & \textbf{8} & \textbf{18} & \textbf{41\%}\\
      \cline{2-16}
      \noalign{\vskip-2\tabcolsep \vskip-3\arrayrulewidth \vskip\doublerulesep} \\
      \cline{2-16}
      & \textbf{RC2\_2\_1} & 372 & 32 & 6 & 38
                           & 586 (378) & 18 & 10 & 28 & \textbf{28\%}
                           & 568 (376) & 17 & 9 & 26 & \textbf{32\%}\\
      & \textbf{RC2\_2\_2} & 218 & 17 & 3 & 20
                           & 332 (179) & 12 & 5 & 17 & \textbf{15\%}
                           & 282 (142) & 11 & 5 & 15 & \textbf{25\%}\\
      & \textbf{RC2\_2\_3} & 852 & 40 & 11 & 51
                           & 912 (695) & 28 & 11 & 39 & \textbf{23\%}
                           & 670 (446) & 22 & 10 & 32 & \textbf{38\%}\\
      & \textbf{RC2\_2\_4} & 609 & 34 & 13 & 47
                           & 675 (476) & 18 & 14 & 33 & \textbf{30\%}
                           & 531 (352) & 17 & 13 & 29 & \textbf{37\%}\\
      & \textbf{RC2\_2\_5} & 313 & 32 & 6 & 39
                           & 451 (226) & 15 & 7 & 23 & \textbf{41\%}
                           & 497 (270) & 20 & 9 & 29 & \textbf{25\%}\\
      \cline{2-16}
      & \textbf{Average} 
      & \textbf{473} & \textbf{31} & \textbf{8} & \textbf{39} 
      & \textbf{591 (391)} & \textbf{18} & \textbf{10} & \textbf{28} & \textbf{28\%} 
      & \textbf{510 (317)} & \textbf{17} & \textbf{9} & \textbf{26} & \textbf{31\%}\\
      \cline{2-16}
      \noalign{\vskip-2\tabcolsep \vskip-3\arrayrulewidth \vskip\doublerulesep} \\
      \cline{2-16}
      & \textbf{C2\_2\_1} & 760 & 72 & 11 & 83
                           & 924 (552) & 54 & 15 & 69 & \textbf{17\%}
                           & 1020 (676) & 62 & 15 & 77 & \textbf{6\%}\\
      & \textbf{C2\_2\_2} & 316 & 121 & 6 & 126
                           & 416 (242) & 76 & 8 & 84 & \textbf{34\%}
                           & 385 (210) & 71 & 8 & 79 & \textbf{38\%}\\
      & \textbf{C2\_2\_3} & 272 & 54 & 7 & 61
                           & 505 (208) & 49 & 11 & 60 & \textbf{2\%}
                           & 489 (190) & 51 & 11 & 62 & \textbf{-1\%}\\
      & \textbf{C2\_2\_4} & 181 & 65 & 5 & 70
                           & 442 (113) & 61 & 11 & 72 & \textbf{-3\%}
                           & 464 (122) & 65 & 12 & 77 & \textbf{-10\%}\\
      & \textbf{C2\_2\_5} & 237 & 140 & 6 & 146
                           & 501 (159) & 127 & 12 & 139 & \textbf{5\%}
                           & 479 (174) & 121 & 11 & 132 & \textbf{10\%}\\
      \cline{2-16}
      & \textbf{Average} 
      & \textbf{353} & \textbf{90} & \textbf{7} & \textbf{97} 
      & \textbf{558 (255)} & \textbf{73} & \textbf{11} & \textbf{84} & \textbf{11\%} 
      & \textbf{567 (274)} & \textbf{74} & \textbf{11} & \textbf{85} & \textbf{10\%}\\
      \cline{1-16}
      \noalign{\vskip-2\tabcolsep \vskip-3\arrayrulewidth \vskip\doublerulesep} \\
      \cline{1-16}
    \multirow{28}{*}{\rotatebox[origin=c]{90}{\textbf{400 clients}}} & \textbf{R2\_4\_1} & 373 & 66 & 119 & 185
                           & 647 (180) & 39 & 119 & 158 & \textbf{15\%}
                           & 678 (230) & 41 & 122 & 163 & \textbf{12\%}\\
      & \textbf{R2\_4\_2} & 289 & 45 & 84 & 129
                           & 436 (127) & 24 & 81 & 105 & \textbf{19\%}
                           & 429 (117) & 24 & 81 & 105 & \textbf{19\%}\\
      & \textbf{R2\_4\_3} & 235 & 97 & 37 & 134
                           & 392 (111) & 40 & 54 & 94 & \textbf{30\%}
                           & 440 (157) & 44 & 55 & 99 & \textbf{0\%}\\
      & \textbf{R2\_4\_4} & 239 & 142 & 20 & 162
                           & 297 (106) & 40 & 29 & 69 & \textbf{57\%}
                           & 314 (134) & 48 & 29 & 77 & \textbf{52\%}\\
      & \textbf{R2\_4\_5} & 361 & 51 & 127 & 178
                           & 619 (165) & 32 & 114 & 146 & \textbf{18\%}
                           & 573 (159) & 30 & 114 & 144 & \textbf{19\%}\\
      & \textbf{R2\_4\_6} & 277 & 44 & 78 & 122
                           & 523 (196) & 32 & 94 & 126 & \textbf{-3\%}
                           & 460 (139) & 25 & 81 & 106 & \textbf{13\%}\\
      & \textbf{R2\_4\_7} & 218 & 81 & 44 & 125
                       & 390 (114) & 38 & 60 & 98 & \textbf{22\%}
                       & 384 (116) & 31 & 53 & 84 & \textbf{33\%}\\
      & \textbf{R2\_4\_8} & 242 & 144 & 23 & 167
                       & 328 (131) & 60 & 33 & 93 & \textbf{44\%}
                       & 309 (119) & 45 & 31 & 76 & \textbf{54\%}\\
      & \textbf{R2\_4\_9} & 322 & 101 & 116 & 217
                       & 615 (168) & 56 & 131 & 187 & \textbf{14\%}
                       & 632 (188) & 50 & 121 & 171 & \textbf{21\%}\\
      & \textbf{R2\_4\_10} & 380 & 154 & 108 & 262
                       & 691 (230) & 76 & 119 & 195 & \textbf{26\%}
                       & 719 (201) & 64 & 116 & 180 & \textbf{31\%}\\
      \cline{2-16}
      & \textbf{Average} 
      & \textbf{294} & \textbf{93} & \textbf{76} & \textbf{168} 
      & \textbf{494 (153)} & \textbf{44} & \textbf{83} & \textbf{127} & \textbf{24\%} 
      & \textbf{494 (156)} & \textbf{40} & \textbf{80} & \textbf{121} & \textbf{28\%}\\
      \cline{2-16}
      \noalign{\vskip-2\tabcolsep \vskip-3\arrayrulewidth \vskip\doublerulesep} \\
      \cline{2-16}
    & \textbf{RC2\_4\_1} & 390 & 45 & 102 & 147
                           & 762 (232) & 42 & 124 & 166 & \textbf{-13\%}
                           & 801 (310) & 50 & 132 & 182 & \textbf{-24\%}\\
      & \textbf{RC2\_4\_2} & 295 & 37 & 67 & 104
                           & 630 (210) & 37 & 76 & 113 & \textbf{-9\%}
                           & 651 (309) & 43 & 89 & 132 & \textbf{-27\%}\\
      & \textbf{RC2\_4\_3} & 301 & 108 & 27 & 135
                           & 446 (188) & 52 & 35 & 87 & \textbf{36\%}
                           & 486 (206) & 53 & 36 & 89 & \textbf{34\%}\\
      & \textbf{RC2\_4\_4} & 183 & 81 & 16 & 97
                           & 263 (76) & 29 & 23 & 52 & \textbf{46\%}
                           & 278 (86) & 31 & 23 & 54 & \textbf{44\%}\\
      & \textbf{RC2\_4\_5} & 1164 & 84 & 186 & 270
                           & 1600 (882) & 83 & 196 & 279 & \textbf{-3\%}
                           & 1810 (1084) & 101 & 251 & 352 & \textbf{-30\%}\\
      & \textbf{RC2\_4\_6} & 434 & 98 & 118 & 216
                           & 773 (244) & 62 & 130 & 192 & \textbf{11\%}
                           & 867 (314) & 62 & 122 & 184 & \textbf{15\%}\\
      & \textbf{RC2\_4\_7} & 415 & 100 & 98 & 198
                       & 792 (348) & 79 & 118 & 197 & \textbf{1\%}
                       & 674 (274) & 60 & 119 & 179 & \textbf{10\%}\\
      & \textbf{RC2\_4\_8} & 435 & 90 & 72 & 162
                       & 834 (357) & 78 & 84 & 162 & \textbf{0\%}
                       & 669 (240) & 56 & 74 & 130 & \textbf{20\%}\\
      & \textbf{RC2\_4\_9} & 371 & 83 & 91 & 174
                       & 723 (304) & 67 & 116 & 183 & \textbf{-5\%}
                       & 952 (561) & 111 & 108 & 219 & \textbf{-26\%}\\
      & \textbf{RC2\_4\_10} & 438 & 40 & 91 & 131
                       & 679 (274) & 39 & 118 & 157 & \textbf{-20\%}
                       & 689 (266) & 35 & 91 & 126 & \textbf{4\%}\\
      \cline{2-16}
      & \textbf{Average} 
      & \textbf{443} & \textbf{77} & \textbf{87} & \textbf{163} 
      & \textbf{750 (312)} & \textbf{57} & \textbf{102} & \textbf{159} & \textbf{4\%} 
      & \textbf{788 (365)} & \textbf{60} & \textbf{105} & \textbf{165} & \textbf{2\%}\\
      \cline{2-16}
      \noalign{\vskip-2\tabcolsep \vskip-3\arrayrulewidth \vskip\doublerulesep} \\
      \cline{2-16}
    & \textbf{C2\_4\_1} & 966 & 144 & 183 & 327
                           & 1692 (669) & 162 & 190 & 352 & \textbf{-8\%}
                           & 1785 (768) & 170 & 186 & 356 & \textbf{-9\%}\\
      & \textbf{C2\_4\_2} & 445 & 71 & 88 & 159
                           & 867 (438) & 99 & 117 & 216 & \textbf{-36\%}
                           & 726 (239) & 68 & 110 & 178 & \textbf{-12\%}\\
      & \textbf{C2\_4\_3} & 484 & 176 & 55 & 231
                           & 741 (333) & 215 & 78 & 293 & \textbf{-27\%}
                           & 913 (463) & 210 & 87 & 297 & \textbf{-29\%}\\
      & \textbf{C2\_4\_4} & 332 & 267 & 30 & 297
                           & 522 (249) & 327 & 44 & 371 & \textbf{-25\%}
                           & 536 (287) & 290 & 47 & 337 & \textbf{-13\%}\\
      & \textbf{C2\_4\_5} & 417 & 194 & 63 & 257
                           & 887 (266) & 248 & 101 & 349 & \textbf{-36\%}
                           & 895 (278) & 222 & 101 & 323 & \textbf{-26\%}\\
      \cline{2-16}
      & \textbf{Average} 
      & \textbf{529} & \textbf{170} & \textbf{84} & \textbf{254} 
      & \textbf{942 (391)} & \textbf{210} & \textbf{106} & \textbf{316} & \textbf{-26\%} 
      & \textbf{971 (407)} & \textbf{192} & \textbf{106} & \textbf{298} & \textbf{-17\%}\\
      \hline
    \end{tabular}}
    \caption{Additional VRPTW results.}
    \label{tab:table-results-vrptw-2}
  \end{center}
\end{table}

From the results, we can see that the gains obtained by the two new strategies are quite comparable. As before, a decrease in performance is noticed when dealing with the $RC2$ and $C2$ instances. On average, ML-RedCost-S-2 seems to be slightly better as it gives a few larger gains for some groups of instances. Overall, the results show that it is possible to get improvements on top of RedCost-S, up to $41\%$ for the 200-client $R2$ instances and $28\%$ for the 400-client ones, especially when the PP is much more time consuming than the RMP. However, it seems to be less effective on $RC2$ and $C2$ instances.

\subsubsection{\label{subsubsec:case2_additional_exp}Additional experiments}

Another attempt to deal with the inefficiency of the model on the $C2$ instances was conducted. On the ML side, one can think that the $C2$ instances are quite different from the other groups (i.e., $RC$ and $R$), and that despite the overall accuracy of $90\%$ obtained, the accuracy is perhaps less for the $C2$ instances. It turns out that this is not the case, the metrics reported for the C2 instances alone are the following: $94\%$ recall, $79\%$ TNR and $87\%$ accuracy. It is true that the accuracy is 3\% below the average, but these are still solid results. Moreover, even with a new model that was only trained on the C2 instances, the results obtained when integrated to the CG algorithm were not much better.

\section{\label{sec:conclusion}Conclusion}

In this paper, a new pricing heuristic based on ML was presented. The goal is to speed up the CG method for problems where the PP is a SPPRC or one of its variants defined on a network, and solved using a labeling algorithm. The method consists in reducing the network size, by keeping only the most promising arcs that have a high chance to be part of good columns. The ML model is trained on the data collected from previous executions through a supervised learning approach. 
The arcs selected by the trained model are used to build a new reduced network that is on average $15\%$ to $25\%$ the size of the original network. The new network is used at each iteration as long as it generates a satisfactory number of columns. If not, the full network is used instead, especially in the last iterations. 

The approach is fairly general and can be used for different problems. We chose to demonstrate it on two well-known problems, the VCSP and VRPTW. 
For the VCSP, the selection was limited to the arcs representing walking movements, since they represent more than 95\% of the arcs in the network. Whereas for the VRPTW, all the arcs were targeted by the selection. The results showed reductions in computing time of up to 27\% for the VCSP, and 40\% for the VRPTW. The resulting ML models have also shown the ability to generalize to new instances not encountered during the training phase. However, the ML model had some limitations when dealing with the $C$ instances of the VRPTW.

\vspace*{5mm}

\textbf{Acknowledgement:} We thank Giro Inc.~and the Natural Sciences and Engineering Research Council of Canada for their financial support through the grant CRDPJ 520349-17. 

\bibliography{main} 

\begin{thebibliography}{26}
\providecommand{\natexlab}[1]{#1}
\providecommand{\url}[1]{\texttt{#1}}
\expandafter\ifx\csname urlstyle\endcsname\relax
  \providecommand{\doi}[1]{doi: #1}\else
  \providecommand{\doi}{doi: \begingroup \urlstyle{rm}\Url}\fi

\bibitem[Alvarez et~al.(2017)Alvarez, Louveaux, and Wehenkel]{alvarez}
Alejandro Alvarez, Quentin Louveaux, and Louis Wehenkel.
\newblock A machine learning-based approximation of strong branching.
\newblock \emph{INFORMS Journal on Computing}, 29:\penalty0 185--195, 01 2017.
\newblock \doi{10.1287/ijoc.2016.0723}.

\bibitem[Baldacci et~al.(2011)Baldacci, Mingozzi, and
  Roberti]{baldacci-ng-routes-2011}
Roberto Baldacci, Aristide Mingozzi, and Roberto Roberti.
\newblock New route relaxation and pricing strategies for the vehicle routing
  problem.
\newblock \emph{Operations Research}, 59\penalty0 (5):\penalty0 1269--1283,
  2011.
\newblock \doi{10.1287/opre.1110.0975}.

\bibitem[Barnhart et~al.(1996)Barnhart, Johnson, Nemhauser, Savelsbergh, and
  Vance]{barnhart1998}
Cynthia Barnhart, Ellis~L. Johnson, George~L. Nemhauser, Martin W.~P.
  Savelsbergh, and Pamela~H. Vance.
\newblock Branch-and-price: Column generation for solving huge integer
  programs.
\newblock \emph{Operations Research}, 46:\penalty0 316--329, 1996.

\bibitem[Bengio et~al.(2020)Bengio, Lodi, and Prouvost]{bengio2018}
Yoshua Bengio, Andrea Lodi, and Antoine Prouvost.
\newblock Machine learning for combinatorial optimization: a methodological
  tour d'horizon.
\newblock \emph{European Journal of Operational Research}, 2020.
\newblock \doi{10.1016/j.ejor.2020.07.063}.

\bibitem[Contardo et~al.(2015)Contardo, Desaulniers, and Lessard]{Contardo2015}
Claudio Contardo, Guy Desaulniers, and Fran{\c{c}}ois Lessard.
\newblock Reaching the elementary lower bound in the vehicle routing problem
  with time windows.
\newblock \emph{Networks}, 65\penalty0 (1):\penalty0 88--99, 2015.

\bibitem[Costa et~al.(2019)Costa, Contardo, and
  Desaulniers]{costa-contardo-desaulniers-2019}
Luciano Costa, Claudio Contardo, and Guy Desaulniers.
\newblock Exact branch-price-and-cut algorithms for vehicle routing.
\newblock \emph{Transportation Science}, 53\penalty0 (4):\penalty0 946--985,
  2019.
\newblock \doi{10.1287/trsc.2018.0878}.

\bibitem[Desaulniers et~al.(1998)Desaulniers, Desrosiers, loachim, Solomon,
  Soumis, and Villeneuve]{Desaulniers1998}
Guy Desaulniers, Jacques Desrosiers, Irina loachim, Marius~M. Solomon,
  Fran{\c{c}}ois Soumis, and Daniel Villeneuve.
\newblock \emph{A Unified Framework for Deterministic Time Constrained Vehicle
  Routing and Crew Scheduling Problems}, pages 57--93.
\newblock Springer US, Boston, MA, 1998.
\newblock ISBN 978-1-4615-5755-5.
\newblock \doi{10.1007/978-1-4615-5755-5_3}.
\newblock URL \url{https://doi.org/10.1007/978-1-4615-5755-5_3}.

\bibitem[Desaulniers et~al.(2005)Desaulniers, Desrosiers, and
  Solomon]{desaulniers2005a}
Guy Desaulniers, Jacques Desrosiers, and Marius~M. Solomon.
\newblock \emph{Column generation}.
\newblock Springer, New York, January 2005.

\bibitem[Desaulniers et~al.(2008)Desaulniers, Lessard, and
  Hadjar]{desaulniers-lessard-tabu-kpath}
Guy Desaulniers, François Lessard, and Ahmed Hadjar.
\newblock Tabu search, partial elementarity, and generalized k-path
  inequalities for the vehicle routing problem with time windows.
\newblock \emph{Transportation Science}, 42\penalty0 (3):\penalty0 387--404,
  2008.
\newblock \doi{10.1287/trsc.1070.0223}.

\bibitem[Desrochers and Soumis(1988)]{desrochers-spprc-1988}
M.~Desrochers and F.~Soumis.
\newblock A generalized permanent labeling algorithm for the shortest path
  problem with time windows.
\newblock \emph{Information Systems and Operations Research}, 26\penalty0
  (3):\penalty0 191--212, 1988.

\bibitem[Desrochers et~al.(1992)Desrochers, Desrosiers, and
  Solomon]{desrochers-vrptw-1992}
Martin Desrochers, Jacques Desrosiers, and Marius Solomon.
\newblock A new optimization algorithm for the vehicle routing problem with
  time windows.
\newblock \emph{Operations Research}, 40\penalty0 (2):\penalty0 342--354, 1992.
\newblock \doi{10.1287/opre.40.2.342}.

\bibitem[Feillet et~al.(2004)Feillet, Dejax, Gendreau, and
  Gueguen]{Feillet2004}
Dominique Feillet, Pierre Dejax, Michel Gendreau, and Cyrille Gueguen.
\newblock An exact algorithm for the elementary shortest path problem with
  resource constraints: Application to some vehicle routing problems.
\newblock \emph{Networks}, 44\penalty0 (3):\penalty0 216--229, 2004.

\bibitem[Freling et~al.(1999)Freling, Wagelmans, and Paix{\~a}o]{freling1999}
Richard Freling, Albert P.~M. Wagelmans, and Jos{\'e} M.~Pinto Paix{\~a}o.
\newblock An overview of models and techniques for integrating vehicle and crew
  scheduling.
\newblock In Nigel H.~M. Wilson, editor, \emph{Computer-Aided Transit
  Scheduling}, pages 441--460, Berlin, Heidelberg, 1999. Springer Berlin
  Heidelberg.
\newblock ISBN 978-3-642-85970-0.

\bibitem[Fukasawa et~al.(2006)Fukasawa, Lysgaard, Poggi~de Arag{\~a}o, Reis,
  Uchoa, and Werneck]{fukasawa2006}
Ricardo Fukasawa, Jens Lysgaard, Marcus Poggi~de Arag{\~a}o, Marcelo Reis,
  Eduardo Uchoa, and Renato~F. Werneck.
\newblock Robust branch-and-cut-and-price for the capacitated vehicle routing
  problem.
\newblock \emph{Mathematical Programming}, 106:\penalty0 491–511, 2006.
\newblock \doi{10.1007/s10107-005-0644-x}.

\bibitem[Gamache et~al.(1999)Gamache, Soumis, Marquis, and
  Desrosiers]{gamache-rostering-99}
Michel Gamache, François Soumis, Gérald Marquis, and Jacques Desrosiers.
\newblock A column generation approach for large-scale aircrew rostering
  problems.
\newblock \emph{Operations Research}, 47\penalty0 (2):\penalty0 247--263, 1999.
\newblock \doi{10.1287/opre.47.2.247}.

\bibitem[Gasse et~al.(2019)Gasse, Chételat, Ferroni, Charlin, and
  Lodi]{gasse2019}
Maxime Gasse, Didier Chételat, Nicola Ferroni, Laurent Charlin, and Andrea
  Lodi.
\newblock Exact combinatorial optimization with graph convolutional neural
  networks.
\newblock In \emph{Advances in Neural Information Processing Systems 32}, 2019.

\bibitem[Haase et~al.(2001)Haase, Desaulniers, and
  Desrosiers]{haase-desaulniers-desrosiers}
Knut Haase, Guy Desaulniers, and Jacques Desrosiers.
\newblock Simultaneous vehicle and crew scheduling in urban mass transit
  systems.
\newblock \emph{Transportation Science}, 35\penalty0 (3):\penalty0 286--303,
  2001.
\newblock \doi{10.1287/trsc.35.3.286.10153}.

\bibitem[Homberger and Gehring(2005)]{homberger-2005}
Jörg Homberger and Hermann Gehring.
\newblock A two-phase hybrid metaheuristic for the vehicle routing problem with
  time windows.
\newblock \emph{European Journal of Operational Research}, 162\penalty0
  (1):\penalty0 220 -- 238, 2005.
\newblock \doi{https://doi.org/10.1016/j.ejor.2004.01.027}.

\bibitem[Irnich(2008)]{irnich2008}
Stefan Irnich.
\newblock Resource extension functions: properties, inversion, and
  generalization to segments.
\newblock \emph{OR Spectrum}, 30\penalty0 (1):\penalty0 113--148, Jan 2008.
\newblock \doi{10.1007/s00291-007-0083-6}.

\bibitem[Irnich and Desaulniers(2005)]{irnich-spprc}
Stefan Irnich and Guy Desaulniers.
\newblock Shortest path problems with resource constraints.
\newblock In Guy Desaulniers, Jacques Desrosiers, and Marius~M. Solomon,
  editors, \emph{Column Generation}, pages 33--65. Springer US, Boston, MA,
  2005.

\bibitem[Irnich and Villeneuve(2006)]{irnich-spprc-k-cyc}
Stefan Irnich and Daniel Villeneuve.
\newblock The shortest-path problem with resource constraints and k-cycle
  elimination for k >= 3.
\newblock \emph{INFORMS Journal on Computing}, 18\penalty0 (3):\penalty0
  391--406, 2006.
\newblock \doi{10.1287/ijoc.1040.0117}.

\bibitem[Khalil et~al.(2016)Khalil, Bodic, Song, Nemhauser, and
  Dilkina]{Khalil2016}
Elias~B. Khalil, Pierre~Le Bodic, Le~Song, George Nemhauser, and Bistra
  Dilkina.
\newblock Learning to branch in mixed integer programming.
\newblock In \emph{Proceedings of the Thirtieth AAAI Conference on Artificial
  Intelligence}, AAAI'16, pages 724--731. AAAI Press, 2016.

\bibitem[Lodi and Zarpellon(2017)]{lodi2017}
Andrea Lodi and Giulia Zarpellon.
\newblock On learning and branching: a survey.
\newblock \emph{TOP}, 25\penalty0 (2):\penalty0 207--236, Jul 2017.
\newblock \doi{10.1007/s11750-017-0451-6}.

\bibitem[Morabit et~al.(2021)Morabit, Desaulniers, and
  Lodi]{mouad-columnselection}
Mouad Morabit, Guy Desaulniers, and Andrea Lodi.
\newblock Machine-learning–based column selection for column generation.
\newblock \emph{Transportation Science}, 55\penalty0 (4):\penalty0 815--831,
  2021.
\newblock \doi{10.1287/trsc.2021.1045}.

\bibitem[Righini and Salani(2006)]{righini-bidir-espprc}
Giovanni Righini and Matteo Salani.
\newblock Symmetry helps: Bounded bi-directional dynamic programming for the
  elementary shortest path problem with resource constraints.
\newblock \emph{Discrete Optimization}, 3\penalty0 (3):\penalty0 255--273,
  2006.
\newblock \doi{https://doi.org/10.1016/j.disopt.2006.05.007}.
\newblock Graphs and Combinatorial Optimization.

\bibitem[Václavík et~al.(2018)Václavík, Novák, Šůcha, and
  Hanzálek]{vaclavik2018}
Roman Václavík, Antonín Novák, Přemysl Šůcha, and Zdeněk Hanzálek.
\newblock Accelerating the branch-and-price algorithm using machine learning.
\newblock \emph{European Journal of Operational Research}, 271\penalty0
  (3):\penalty0 1055 -- 1069, 2018.
\newblock \doi{https://doi.org/10.1016/j.ejor.2018.05.046}.

\end{thebibliography}
\end{document}